\documentclass{article}
\usepackage[a4paper, total={6.5in, 9in}]{geometry}
\usepackage[table]{xcolor}

\usepackage[scr=rsfs]{mathalpha}

\usepackage{tikz}
\usetikzlibrary{}
\usetikzlibrary{matrix,chains,positioning,decorations.pathreplacing,arrows}
\usetikzlibrary{positioning,calc}
\usetikzlibrary{graphs,graphs.standard,quotes}

\usepackage{pgfplots}
\pgfplotsset{compat=1.17} 
\usepackage[utf8]{inputenc}
\usepackage[square,sort,comma,numbers]{natbib}
\usepackage{graphicx}
\graphicspath{ {./images/} }
\usepackage{amsmath}
\usepackage{amsfonts}
\usepackage{algorithm}
\usepackage{algorithmic}
\usepackage{mathtools}
\usepackage{mathrsfs}
\usepackage{tabularx}
\usepackage{multirow}
\usepackage{float}
\usepackage{subcaption}
\usepackage[font={footnotesize}]{caption}
\usepackage{stackengine}
\usepackage{url}
\usepackage{booktabs}
\usepackage{makecell}
\usepackage{lineno}
\usepackage{enumitem}
\usepackage{array}
\usepackage[hidelinks]{hyperref}

\newlength{\bibitemsep}
\newlength{\bibparskip}
\let\oldthebibliography\thebibliography
\renewcommand\thebibliography[1]{%
  \footnotesize
  \oldthebibliography{#1}%
  \setlength{\parskip}{\bibitemsep}%
  \setlength{\itemsep}{\bibparskip}%
}

\newtheorem{definition}{Definition}[section]
\newtheorem{remark}{Remark}

\newcommand{\norm}[1]{\left\|#1\right\|}

\DeclareMathOperator*{\argmin}{argmin}

\title{A Deep Learning algorithm to accelerate Algebraic Multigrid methods in Finite Element solvers of 3D PDEs}
\author{Matteo Caldana$^{a,}$\thanks{Corresponding author: {\tt matteo.caldana@polimi.it}},
Paola F. Antonietti$^a$, 
Luca Dede'$^a$\\[0.3cm]
\small\textit{$^a$MOX, Dipartimento di Matematica, 
Politecnico di Milano,
Piazza Leonardo da Vinci 32,
20133 Milano, Italy
}}

{\date{\small\today}}

\begin{document}

\maketitle
{{\centering\subsection*{Abstract}}
Algebraic multigrid (AMG) methods are among the most efficient solvers for linear systems of equations and they are widely used for the solution of problems stemming from the discretization of Partial Differential Equations (PDEs). A severe limitation of AMG methods is the dependence on parameters that require to be fine-tuned. In particular, the strong threshold parameter is the most relevant since it stands at the basis of the construction of successively coarser grids needed by the AMG methods. We introduce a novel deep learning algorithm that minimizes the computational cost of the AMG method when used as a finite element solver. We show that our algorithm requires minimal changes to any existing code. The proposed Artificial Neural Network (ANN) tunes the value of the strong threshold parameter by interpreting the sparse matrix of the linear system as a gray scale image and exploiting a pooling operator to transform it into a small multi-channel image. We experimentally prove that the pooling successfully reduces the computational cost of processing a large sparse matrix and preserves the features needed for the regression task at hand. We train the proposed algorithm on a large dataset containing problems with a highly heterogeneous diffusion coefficient defined in different three-dimensional geometries and discretized with unstructured grids and linear elasticity problems with a highly heterogeneous Young's modulus. When tested on problems with coefficients or geometries not present in the training dataset, our approach reduces the computational time by up to 30\%.
}
\vspace{0.4cm}
\paragraph{Key words:} algebraic multigrid methods, partial differential equations, finite element method, elliptic problems, deep learning, convolutional neural networks.
\paragraph{AMS subject classification:} 65N22, 65N30, 65N55, 68T01


\section{Introduction}

Geometric multigrid methods \cite{bramble2019multigrid, wesseling2004introduction} are among the state-of-the-art solvers for large linear systems that come from the discretization of Partial Differential Equations (PDEs). They are applicable to a wide range of discretizations such as polyhedral discontinuous Galerkin \cite{antonietti2017multigrid, antonietti2019v} and virtual elements methods \cite{antonietti2018multigrid, antonietti2023agglomeration}. However, their limitation is that they rely on a sequence of coarser grids that must be available a-priori to solve the problem. The algebraic multigrid (AMG) methods are a highly scalable \cite{baker2012scaling} generalization that is capable of building this hierarchy of grids algebraically. A challenge behind the algebraic construction of a coarser mesh is the selection of the prolongation operator, which in turns relies on seeing the sparse matrix of the linear system as a weighted graph and defining a strong threshold parameter to partition the aforementioned graph. 

The AMG method was first introduced in the 80s \cite{brandt1984algebraic, brandt1986algebraic, brandt1983algebraic, ruge1987algebraic} and gained momentum very quickly is the subsequent years \cite{briggs2000multigrid, falgout2006introduction, trottenberg2000multigrid, vassilevski2008multilevel}. In recent years, modifications have been proposed to the AMG to improve its efficiency, like smoothed aggregation \cite{vanek1996algebraic, van2001convergence} and extend its range of application to different discretization techniques, like discontinuous Galerkin \cite{antonietti2020algebraic}, Ritz-type finite element methods \cite{brezina2001algebraic, chartier2003spectral, jones2001amge} and hybrid discretizations \cite{di2023algebraic, botti2022p}, and to different physics equations like Maxwell's equations \cite{bochev2003improved, kolev2009parallel}, magnetohydrodynamics \cite{adler2016monolithic}, Navier-Stokes equations \cite{weiss1999implicit, raw1996robustness}, linear elasticity \cite{griebel2003algebraic, barnafi2023comparative} and poromechanics \cite{white2019two, arraras2021multigrid}. There is also a wide literature that tackles the AMG from a theoretical point of view, the foundation was laid down in \cite{bramble1991convergence, xu1992iterative, xu2002method, falgout2005two} and the most recent survey is found in \cite{xu2017algebraic}. In particular, we consider the finite element (FE) \cite{quarteroni2008numerical} discretization of diffusion and linear elasticity problems with highly heterogeneous coefficients in three-dimensional space (3D). These problems represent a challenging benchmark since they feature a large number of unknowns and a very ill-conditioned system matrix. Consequently, the efficient application of the AMG method is the key to obtain fast convergence.

In this work, we propose to use Deep Learning (DL) \cite{lecun2015deep} techniques to find the optimal value of the strong threshold parameter, depending on the problem we want to solve. In particular, Artificial Neural Networks (ANNs) have gained widespread popularity for a variety of applications. They are versatile tools that are progressively being used in scientific computing \cite{kutyniok2022theoretical, mishra2018machine, vinuesa2022enhancing}, particularly in the numerical approximation of PDEs and model order reduction \cite{fresca2021comprehensive, regazzoni2019machine}. ANNs can also be utilized within a data-driven framework to provide alternative closure models, which involve learning input-output relationships in complex physical processes \cite{regazzoni2020machine, regazzoni2022machine, berrone2022invariances}. Since we treat the sparse matrix of the linear system like an image, part of the ANN that we employ is comprised by a Convolutional Neural Network (CNN). CNNs are among the most applied neural architectures to analyze visual imagery and achieved breakthrough results \cite{lecun1995convolutional, hu2018squeeze, iandola2016squeezenet, krizhevsky2017imagenet, gu2018recent}. Today, they are only surpassed by visual transformer \cite{dosovitskiy2020image}, which however have a much more complex architecture and are much more difficult to optimize. Moreover, CNNs have already been successfully applied in the field of scientific computing \cite{bhatnagar2019prediction, eichinger2020stationary}.

The combination of finite element/multigrids methods and machine learning has already been used in literature. For instance, in \cite{katrutsa2020black, greenfeld2019learning} there is the first attempt of optimizing the multigrid parameters, in \cite{antonietti2022machine, antonietti2022refinement, antonietti2022agglomeration} DL is used to perform grid refinement and agglomeration, in \cite{luz2020learning, moore2022learning} Graph Neural Networks (GNNs) are used to learn the prolongation operator and in \cite{taghibakhshi2021optimization} reinforcement learning is used to perform graph coarsening, in \cite{heinlein2019machine, heinlein2021combining} DL is used to enhance domain decomposition methods. In this paper, we propose a DL-based algorithm that is able to tune on-the-fight the strong threshold parameter so as to minimize the computational time needed to solve the linear system at hand. The main difference with the previous works is that our algorithm is completely non intrusive: following the approach described in \cite{antonietti2023accelerating}, it does not require any change to existing code (neither to the FE nor the AMG solver). This guarantees a wider range of applicability and means that we can rely on all the classical theoretical results regarding convergence. Namely, an ANN is trained to predict the computational cost of solving a linear system given as input a certain value of the strong threshold parameter and a small multi channel image representing the sparse matrix of the linear system. This representation is built by employing the pooling operator used in CNN. 
However, this algorithm struggles to effectively manage the increased complexity that arises when dealing with 3D problems. In this work we introduce a set of improvements that make the algorithm fast and accurate also for three-dimensional test cases. More precisely, we propose a pre-processing step
that assures data quality, we enhance the architecture of the ANN by adding  an additional output variable measuring the confidence of the neural network on the prediction and we improve the pooling step by means of a four channel tensor. We also show how to speed-up the training by employing layer freezing.
We found that the proposed ANN-enhanced AMG method reduces significantly the computational cost (elapsed time) needed to solve the linear system compared to employing the pre-defined choice of the parameters based on trial-and-error, experience, and literature.

This work is organized as follows. Section~\ref{sec:amg} introduces the mathematical framework of the AMG methods, in particular we define and show the importance of the strong threshold parameter. In Section~\ref{sec:amg-ann} we present our algorithm. We describe the architecture of the ANN and show how to apply the pooling operator to a large sparse matrix. Section~\ref{sec:numerical-results} is devoted to the numerical validation on our method. We first make preliminary sensitivity analysis of the hyperparameters of the ANN. We then apply our algorithm to a family of elliptic problems with a highly heterogeneous diffusion coefficient on structured and unstructured meshes and to a family of linear elasticity problems with a highly heterogeneous Young's modulus. Finally, in Section~\ref{sec:conclusions} we draw our conclusions with future developments.

\section{AMG methods}\label{sec:amg}
To start, we will explain the fundamental concepts and methods that make up AMG. For further information, refer to \cite{ruge1987algebraic, yang2002boomeramg}. Our goal is to solve the linear system $\mathrm{A}\mathbf{u} = \mathbf{f}$ where $\mathrm{A} \in \mathbb{R}^{n_1 \times n_1}$ is a symmetric positive definite (SPD) matrix with entries $a_{ij}$ and $\mathbf{u}, \mathbf{f} \in \mathbb{R}^{n_1}$ are vectors with entries $(\mathbf u)_i$ and $(\mathbf f)_i$, respectively. We define a smoothing scheme as a linear iterative method
\begin{equation}
    \mathbf{u}_{l+1} = \mathbf{u}_{l} + \mathrm{B}(\mathbf{f} - \mathrm A \mathbf{u}_{l}), \quad l \geq 0
    \label{eq:smoother-iteration}
\end{equation}
where $\mathrm B \in \mathbb{R}^{n_1 \times n_1}$ and the initial guess $\mathbf u_0$ are given. We denote $\nu$ applications of (\ref{eq:smoother-iteration}) as
$\texttt{smooth}^\nu(\mathrm A, \mathrm B, \mathbf u_0, \mathbf f)$.
The idea of multigrid methods is that after a certain number of iterations $\nu$ of Eq.~(\ref{eq:smoother-iteration}), the error is more efficiently reduced by projecting the residual equation $\mathrm A \mathbf e = \mathbf r = \mathbf f - \mathrm A \mathbf u_\nu$ on a coarser space, and interpolating the solution back to the original space to apply the correction $\mathbf u_\nu + \mathbf e$. 
If we suppose to have a sequence of interpolation operators $\mathrm I_{k}^{k-1} \in \mathbb R^{n_{k-1} \times n_{k}}$, 
restriction operators $\mathrm I_{k-1}^{k} \in \mathbb R^{n_{k} \times n_{k-1}}$ 
and grid operators $\mathrm A^{(k)} \in \mathbb R^{n_k \times n_k}$ 
for $k = 2, ..., M$ with $\mathrm A^{(1)} = \mathrm A$, 
and of pre- and post- smoothers $\mathrm B_1^{(k)}, \mathrm B_2^{(k)}$ for $k=1,...,M$ with $n_1 > n_2 > ... > n_M$, and a number of pre-smoothing iterations $\nu_1$ and post-smoothing $\nu_2$, the V-cycle multigrid iteration is defined as in Algorithm~\ref{a:v-cycle} (notice the usage of the superscript $^{(k)}$ to indicate the different levels). The AMG method aims at finding the operators needed for this task by just relying upon the original matrix $\mathrm A$. Since $\mathrm A$ is SPD we assume that 
\begin{equation}
    \mathrm{I}_k^{k+1} = (\mathrm I_{k+1}^{k})^\top, \quad \mathrm A^{(k+1)} = \mathrm I_k^{k+1} \mathrm A^{(k)} \mathrm I_{k+1}^{k}, \quad \forall k=1,...,M-1.
\label{eq:op-def-spd}
\end{equation}
Hence, all the operators are defined once we have a recipe to build the interpolation operator.

\begin{algorithm}[t]
\caption{One Iteration of the V-cycle of the AMG method \newline $\mathbf{u}^{(k)} = \texttt{vcycle}^{k}(\mathbf{u}^{(k)}, \mathbf{f}^{(k)}, \{(\mathrm{A}^{(j)}, \mathrm{B}_1^{(j)}, \mathrm{B}_2^{(j)})\}_{j=k}^M, \{(\mathrm{I}^{j+1}_{j}, \mathrm{I}^{j}_{j+1})\}_{j=k}^{M-1}, \nu_1, \nu_2)$}
\label{a:v-cycle}
\begin{algorithmic}[1]
    \IF{$k = M$}
        \STATE $\mathbf{u^{(M)}} = \texttt{gaussian\_elimination}(\mathrm{A}^{(M)}, \mathbf f^{(M)})$\;
    \ELSE
        \STATE $\mathbf{u}^{(k)} \leftarrow \texttt{smooth}^{\nu_1}(\mathrm{A}^{(k)}, \mathrm B_1^{(k)}, \mathbf{u}^{(k)}, \mathbf{f}^{(k)})$
        \STATE $\mathbf{r}^{(k+1)}$ $\leftarrow$ $\mathrm I_k^{k+1} (\mathbf{f}^{(k)} - \mathrm{A}^{(k)} \mathbf{u}^{(k)})$   
        \STATE $\mathbf{e}^{(k+1)}$ $\leftarrow$ $\texttt{vcycle}^{k+1}($\parbox[t]{.6\linewidth}{${\mathbf u}^{(k)}, \mathbf{f}^{(k)}, $ $ \{(\mathrm{A}^{(j)},\mathrm{B}_1^{(j)},\mathrm{B}_2^{(j)})\}_{j=k+1}^M, $ \\ $\{(\mathrm{I}^{j+1}_{j}, \mathrm{I}^{j}_{j+1})\}_{j=k+1}^{M-1}, \nu_1, \nu_2)$}
        \STATE $\mathbf{u}^{(k)}$ $\leftarrow$ $\mathbf{u}^{(k)} + \mathrm{I}_{k+1}^k \mathbf{e}^{(k+1)}$
        \STATE $\mathbf{u}^{(k)} \leftarrow \texttt{smooth}^{\nu_2}(\mathrm{A}^{(k)}, \mathrm B_2^{(k)}, \mathbf{u}^{(k)}, \mathbf{f}^{(k)})$
    \ENDIF
\end{algorithmic}
\end{algorithm}

\subsection{Interpolation operator and coarse-grid selection}\label{sec:interpolation-op}
We consider the interpolation operator between the level $k$ and $k+1$. We assume that we can split the variables into two sets: the one that needs interpolation at the fine level $k$ and the one that are kept at the coarse level $k+1$, namely:
\begin{equation}
(\mathrm{I}_{k+1}^k \mathbf{x})_i = 
\left\{
\begin{matrix}
    (\mathbf{x})_i & \text{ if } i \in \mathcal{C}^k, \\
    \sum_{k \in \mathcal{C}_i^k} \omega_{ij}^k(\mathbf{x})_j & \text{ if } i \in \mathcal{F}^k,
\end{matrix}
\right.
\label{eq:interp-def}
\end{equation}
where $\mathbf x \in \mathbb R^{n_{k+1}}$ is a generic vector, $\mathcal{C}^k / \mathcal{F}^k$ is a coarse/fine partition of the set $\mathcal{N}_k = \{1, ..., n_k\}$, $\mathcal{C}_i^k = \{j \in \mathcal{C}^k: a_{ij} \neq 0\}$ is a subset of $\mathcal{C}^k$ that depends on $i$ and $\omega_{ij}^k$ is a set of weights. 

Since we are working between two levels, from now on, we will omit the superscript $k$. Before being able to define the value of $\omega_{ij}$ we need to define when the unknown $(\mathbf u)_j$ is important in determining the value of $(\mathbf u)_i$. 

\begin{definition}
\label{def:strong_conn}
Given a threshold parameter $0 < \theta \leq 1$, the set of variables on which the variable $i$ strongly depends on is
\begin{linenomath}\begin{equation*}
\mathcal{S}_i = \{j \neq i: - a_{ij} \geq \theta \, \max_{l \neq i} \, \{ - a_{il} \}, j=1,...,n_k \}.
\end{equation*}\end{linenomath}
\end{definition}
We also define the set of variables $j$ that are
influenced by the variable $i$ \begin{linenomath}\begin{equation*}\mathcal{S}_i^\top = \{j: i \in \mathcal{S}_i, j = 1, ..., n_k\}.\end{equation*}\end{linenomath}
The underlying assumption of AMG (see \cite{ruge1987algebraic}) is that the error satisfies
\begin{linenomath}\begin{equation*}\sum_{j\in \mathcal{C}_i} a_{ij} (\mathbf{e})_j+ \sum_{j\in \mathcal{D}_i^s} a_{ij} (\mathbf{e})_j + \sum_{j\in \mathcal{D}_i^w} a_{ij} (\mathbf{e})_j = 0, \quad \forall i = 1, ..., n,\end{equation*}\end{linenomath}
where $\mathcal{D}_i^s = \mathcal{F} \cap \mathcal{S}_i$ and $\mathcal{D}_i^w = \{j \in \mathcal{N}: a_{ij}\neq 0, j \notin \mathcal{S}_i\}$.
This yields the formula for the weights (where $\hat{a}_{ij} = a_{ij}$ if $a_{ij}a_{ii}\leq0$ and zero otherwise)
\begin{equation}
\omega_{ij} = - \frac{1}{a_{ij} + \sum_{l \in \mathcal{D}_i^w}a_{il}} \left(a_{ij}+\sum_{l \in \mathcal{D}_i^s}\frac{a_{il}\hat{a}_{lj}}{\sum_{m\in\mathcal{C}_i}\hat{a}_{lm}}\right)
\label{eq:interp_w_def}
\end{equation}

The last ingredient that we need to build the interpolation operator is the $\mathcal{C}/\mathcal{F}$ splitting. 
Among the several techniques that can be used, we outline the CLJP (Cleary-Luby-Jones-Plassman) algorithm \cite{cleary1998coarse}.  First, we construct the graph of variable dependencies $G=(V, E)$ with vertices $V=\{1, ..., n\}$ and edges $E=\{(i, j)\in V\times V: j \in \mathcal{S}_i\}$. For each vertex, we define the measure $\eta(i) = |\mathcal{S}_i^\top| + \tilde \eta$, where $\tilde \eta$ is a random number in $(0, 1)$ used to break ties. Then, we update $\eta$ and $G$ in the following way until all vertexes are designed as either $\mathcal{C}$ or $\mathcal{F}$. Whenever a update to $\eta(j)$ is such that $\eta(j) < 1$, $j$ is flagged as $\mathcal{F}$.

\begin{enumerate}[label*=\arabic*.,topsep=0pt,itemsep=-1ex,partopsep=1ex,parsep=1ex]
    \item Build the independent set $D=\{i \in V: \eta(i) > \eta(j) \forall j \in \mathcal{S}_i \cap \mathcal{S}_i^\top\}$ for the graph $G$.
    \item For each $i \in D$
    \begin{enumerate}[label*=\arabic*.,topsep=-0.1cm,itemsep=-1ex,partopsep=1ex,parsep=1ex]
    \item For each $j \in \mathcal{S}_i$, decrement $\eta(j)$ and remove the edge $(i, j)$ from $E$. 
    \item For each $j \in \mathcal{S}_i^\top$, remove $(j, i)$ from $E$ 
    \begin{enumerate}[label*=\arabic*.,topsep=-0.1cm,itemsep=-1ex,partopsep=1ex,parsep=1ex]
    \item For each $k \in \mathcal{S}_j^\top \cap \mathcal{S}_i^\top$, decrement $\eta(j)$ and remove $(k, j)$ from $E$. 
    \end{enumerate}
  \end{enumerate}
    \item Every point in $D$ is flagged as $\mathcal C$.
\end{enumerate}
This algorithm is applied at each level until the size of the grid operator $n_k$ is smaller than a certain given value $n_{min}$, that we fix equal to two. Thus, given the strong threshold $\theta$ we are able to construct the succession of operators needed for the V-cycle. Algorithm~\ref{a:amg} showcases the full AMG method.

\begin{algorithm}[t]
    \caption{AMG algorithm \newline $\mathbf{u} = \texttt{AMG}(\mathbf{u}, \mathrm{A}, \mathbf{f}, \theta, \{(\mathrm{B}_1^{(j)}, \mathrm{B}_2^{(j)})\}_{j=k}^M, \nu_1, \nu_2, N_{max}, tol)$}
    \label{a:amg}
    \begin{algorithmic}[1]
        \STATE build $\{\mathcal{C}^k, \mathcal{F}^k\}_{k=1}^M$ using $\theta$ by means of CLJP
    \STATE build the operators $\{\mathrm{I}^{j}_{j+1}\}_{j=k}^{M-1}$ employing Eq.~(\ref{eq:interp-def}) and Eq.~(\ref{eq:interp_w_def})
        \STATE build the operators $\{\mathrm{I}^{j+1}_{j}\}_{j=k}^{M-1}$ and $\{(\mathrm{A}^{(j)}\}_{j=k}^M$ by means of Eq.~(\ref{eq:op-def-spd})
        \WHILE{$k < N_{max} $ \textbf{and} $ \norm{\mathrm{A}\mathbf{u}-\mathbf{f}}/\norm{\mathbf{f}} < tol$}
            \STATE{$\mathbf u \leftarrow  \texttt{vcycle}^{1}(\mathbf{u}, \mathbf{f}, \{(\mathrm{A}^{(j)}, \mathrm{B}_1^{(j)}, \mathrm{B}_2^{(j)})\}_{j=k}^M, \{(\mathrm{I}^{j+1}_{j}, \mathrm{I}^{j}_{j+1})\}_{j=k}^{M-1}, \nu_1, \nu_2)$}
            \STATE $k \leftarrow k + 1$
        \ENDWHILE
    \end{algorithmic}
\end{algorithm}

\section{ANN-enhanced AMG}\label{sec:amg-ann}
We now introduce the ideas behind the proposed AMG-ANN algorithm. We have shown that the strong threshold parameter $\theta$ heavily conditions the construction of the interpolation operator that stands at the basis of the AMG method. In literature, this value is usually fixed to 0.25 and 0.5 when applying the AMG to the solution of linear systems that stem from the discretization of 2D and 3D PDEs, respectively. Our algorithm aims at making an accurate choice of $\theta$ so as to minimize the computational cost of the AMG method. 
Cost can be minimized with respect to various measures, however we choose $\theta$ for two reasons. First, it is the most commonly used parameter, present in almost all the AMG algorithms. Second, from our experiments, it seems that among all parameter we tested, it is the one the most influences the computational cost.
We point out that, in order to keep low the computational costs involved in the optimisation step and to better analyze the proposed algorithm, we decide not to  optimize multiple parameters of the AMG at the same time. Indeed, when optimizing for more than one parameter we would expect to find many local minima and this would hinder the study of the generalization capability of the neural network. Future research efforts are already directed towards adding the ability to choose the best combination among many parameters, such as the type of smoother and the number of smoothing steps.

\begin{figure}[t]
    \makebox[\textwidth][c]{\includegraphics[trim={0 0cm 0 0cm},clip,width=1.0\textwidth]{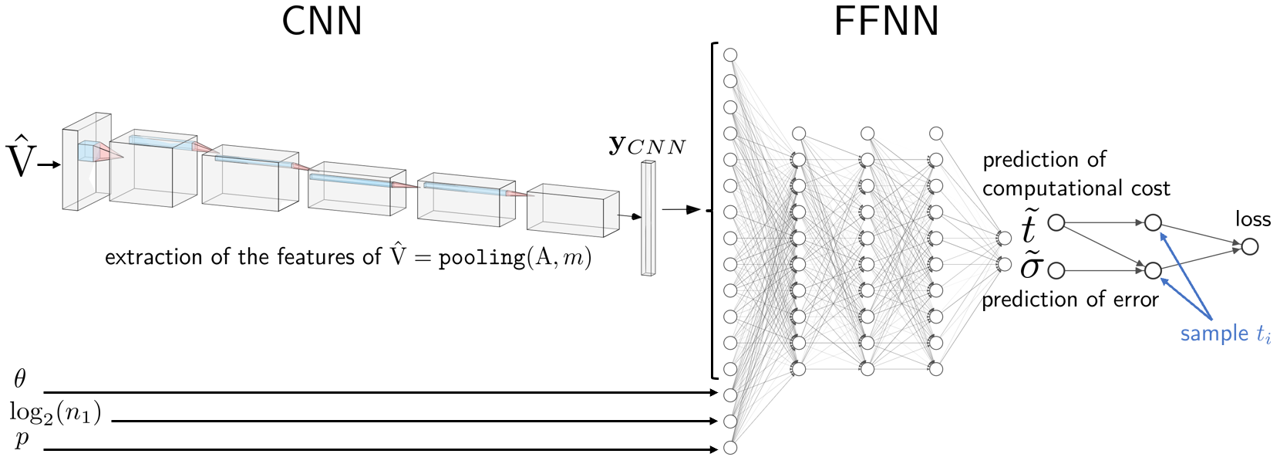}}
    \caption{Architecture of the proposed ANN used to predict the optimal value of the strong threshold parameter $\theta$.}
    \label{fig:architecture}
\end{figure}

To this end, we devise an algorithm that exploits DL techniques to predict the value of $\theta$ that minimizes the elapsed CPU time needed to solve the given linear system $\mathrm A \mathbf u = \mathbf f$. Namely, we propose to use an ANN $\mathscr{F}$ to predict the computational time $t$, that is the number of seconds needed to solve a linear system with the AMG method (including the setup phase), and its confidence of that prediction $\sigma^2$, that is the square of the committed error, given as input a small multi-channel image $\mathrm V$ extracted from the matrix $\mathrm A$, the FE degree $p$, the logarithm of the size of $\mathrm A$ and the value of the strong threshold parameter $\theta$. Then, online, we perform a 1D optimization to find the optimal $\theta^*$ that minimizes the predicted cost $\tilde t$ (we use the tilde to denote quantities approximated by the ANN). The architecture of the ANN is shown in Figure~\ref{fig:architecture}.

This section is structured as follows: in Section~\ref{sec:data-gather}, we briefly review how we collect data for the training of the ANN; in Section~\ref{sec:pooling} and Section~\ref{sec:normalization} we explain how to compress the large sparse matrix $\mathrm A$ so that can be used as input of a  neural network; Section~\ref{sec:ann-architecture} details the architecture of the ANN and how to predict the optimal value of $\theta$. 
Finally, in Section~\ref{sec:amg-ann-algo} we collect together the results of the previous section and describe the AMG-ANN algorithm To make the article self-contained we discuss the basics of DL in Appendix~\ref{apx:dl}.

\subsection{Gathering and smoothing data}\label{sec:data-gather}
\begin{figure}[t]
\setlength{\tabcolsep}{0pt}
    \centering
        \begin{tabular}{cccc}
             \includegraphics[width=0.33\linewidth]{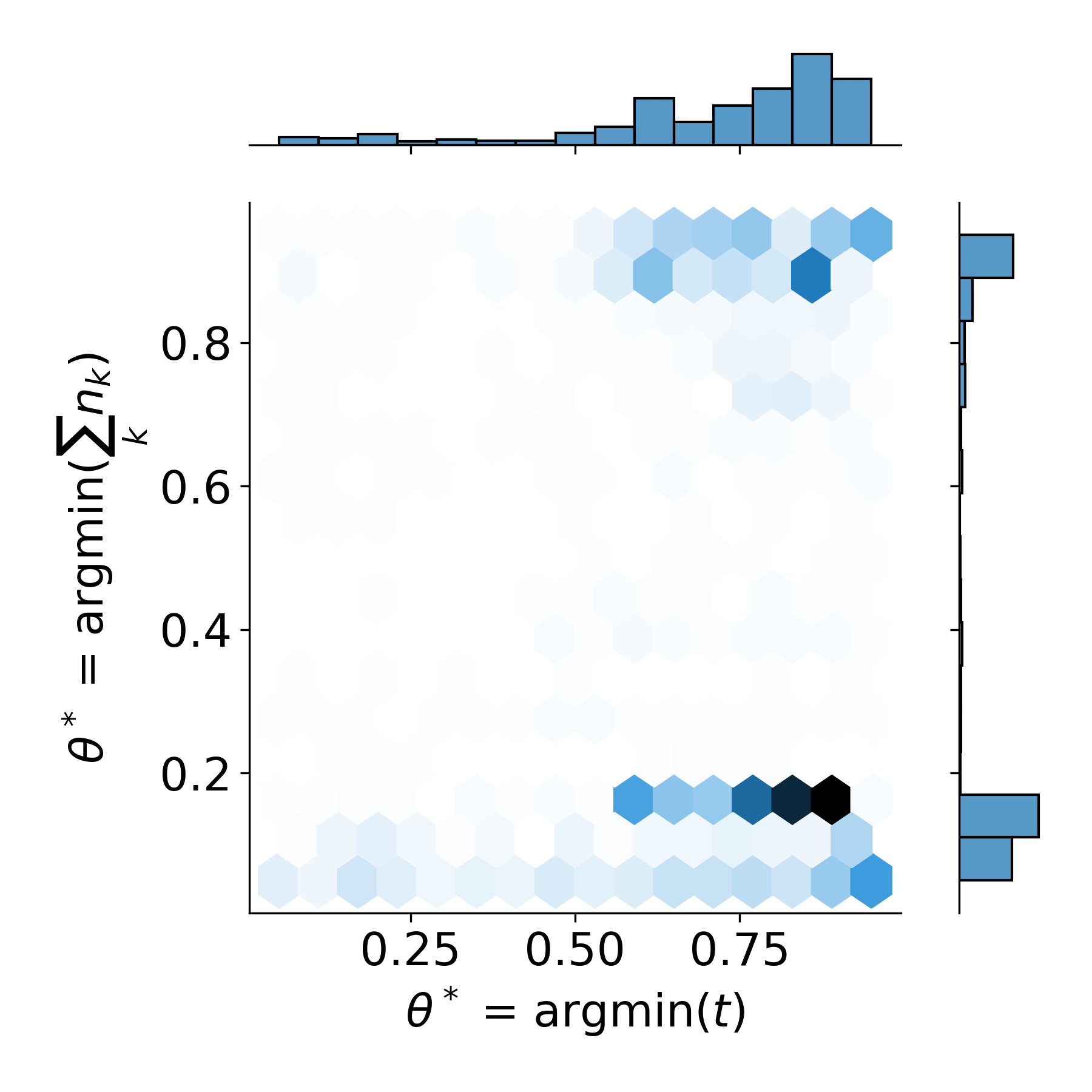}
            & \includegraphics[width=0.33\linewidth]{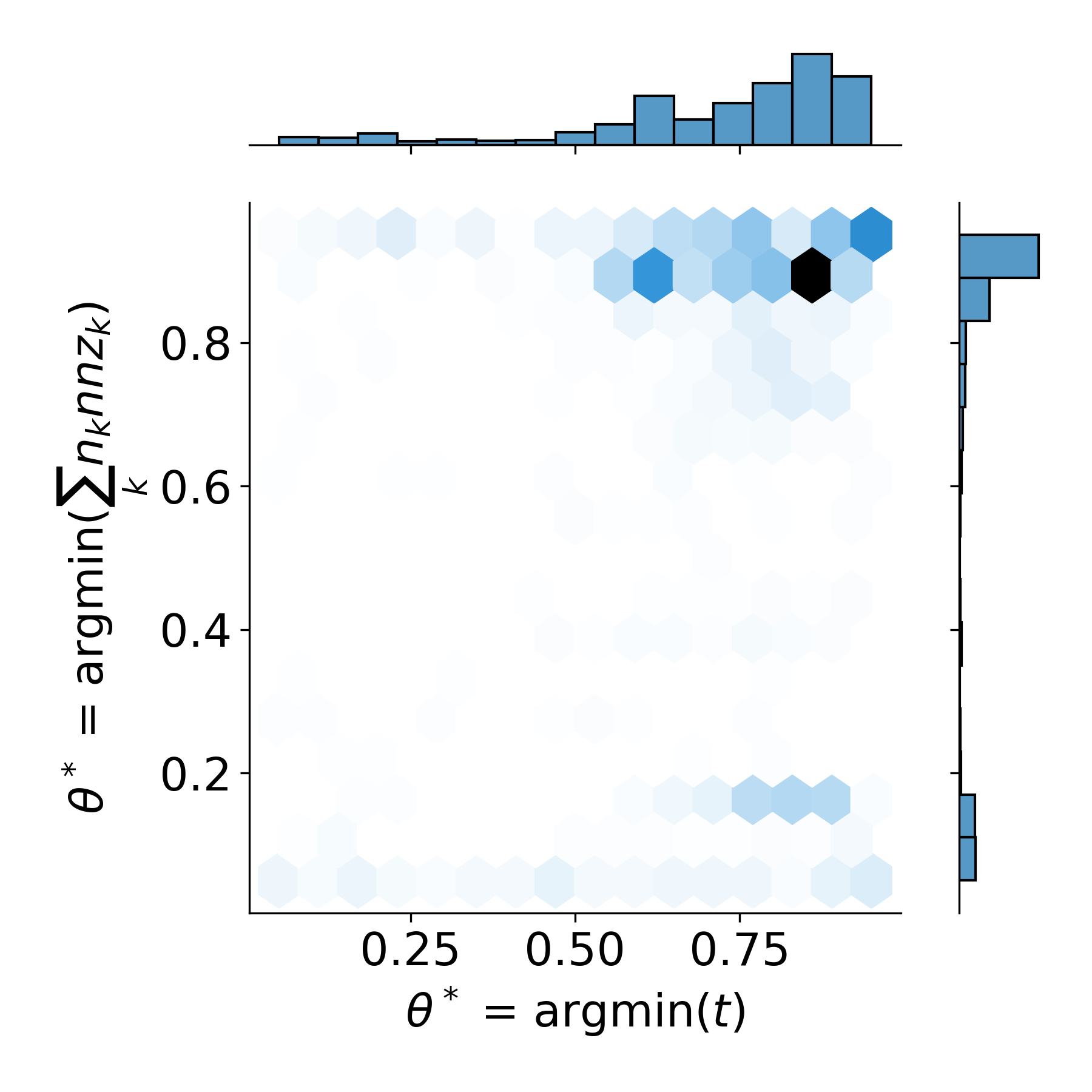}
            & \includegraphics[width=0.33\linewidth]{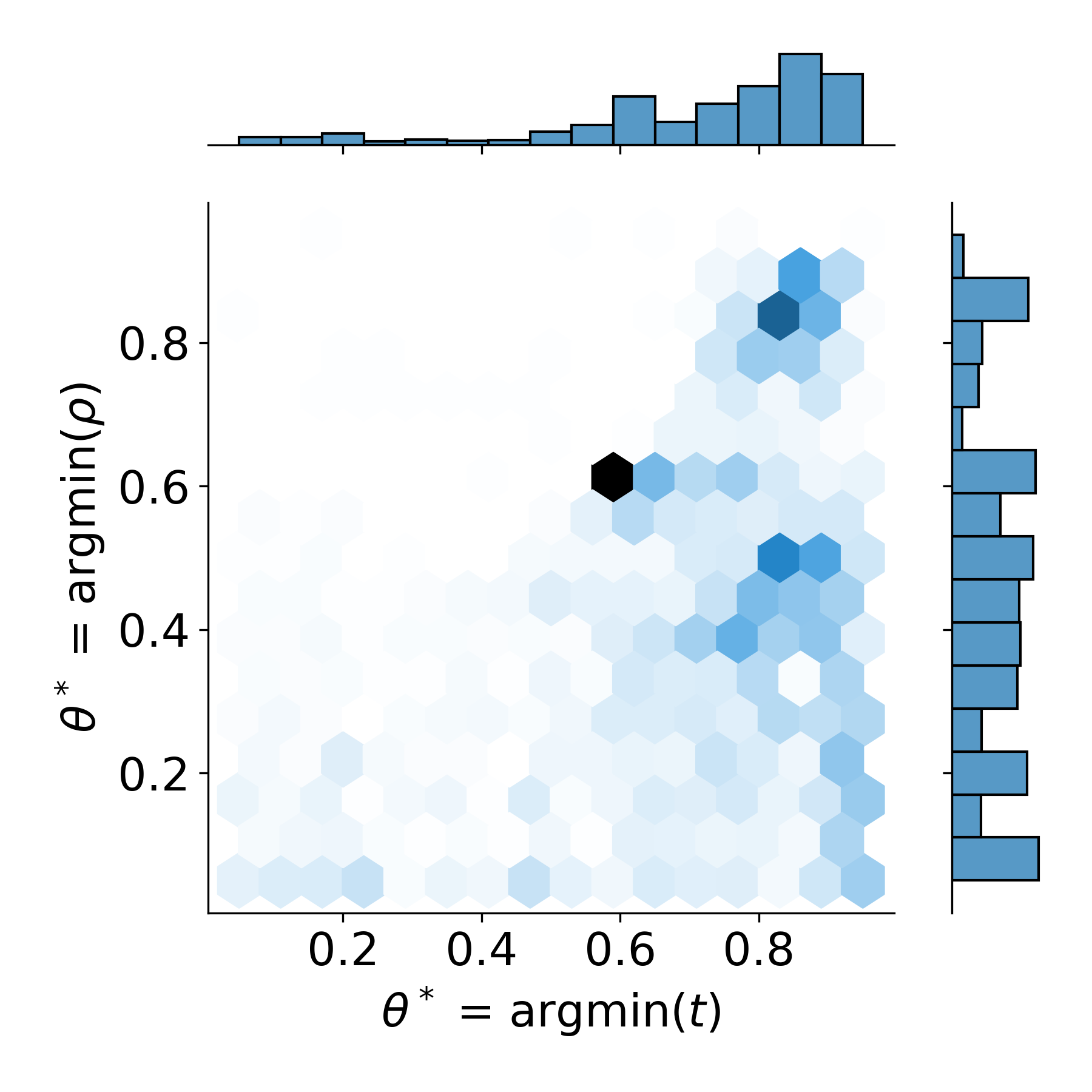}
        \end{tabular}%
    \caption{Joint plot of the optimal strong threshold chosen with respect to different metrics. On the abscissas we have the $\theta^*$ that minimizes the total execution time $t$. None of the performance models considered produces a $\theta^*$ that correlates with it.}\label{fig:argmin-time-vs}
\end{figure}

In this section we describe how we collect the data used to train the neural network. Namely, we are interested in $t$, the sum of the time to setup the AMG method and the execution time to convergence. Since the measurements of $t$ are affected by errors, we need to repeat the measurement. However, this drives up the cost of the method. We would prefer, instead, to use a performance model for predicting the time. We tested several multivariate polynomial models that take into account $n_k$, $nnz_k$ (the size and number of non-zero of $\mathrm A^{(k)}$ at each level $k$, respectively) and the approximate convergence factor $\rho = (\mathbf{r}^{(N_{it})}/\mathbf{r}^{(0)})^{N_{it}}$ (where $N_{it}$ is the number of iterations needed to reach convergence and $\mathbf{r}^{(j)}$ is the residual at the $j$-th iteration), but none were able to consistently give the same choice of $\theta$ that we have by minimizing $t$. We report some relevant tests in a joint plot in Figure~\ref{fig:argmin-time-vs}: we can see that the two quantities do not show any correlation. Our explanation is that there is a complex interaction of factors that drives the total time $t$. Among these factors there is the setup time and the trade-off between performing a larger number of cheaper cycles and a lower number of more expensive cycles.

Then, to minimize the cost of collecting the data we employ the following procedure: we repeat the measurements $r$ times with $r$ between 2 and 100, where $r$ is chosen inversely proportional to the mean elapsed time in the first two measurements. Indeed, the uncertainty in the measurements is caused by the tasks scheduled by the operating system in concurrence with the AMG solver that perturb the CPU load, hence if $t$ is larger the uncertainty is smaller since these perturbations average out. This hypothesis is supported by our measurements, indeed we have observed that the sample variance of our measurements is inversely proportional to $t$ (fixed $r$). The elapsed CPU time $t$ is then defined as the mean of the measurements. 

Moreover, we regularize data using Savitzky-Golay filter \cite{savitzky1964smoothing}. Namely, we fit a polynomial to a successive sub-sets of adjacent data points (called window) by means of least squares. We use a window of size 21 and polynomial degree 7 for uniformly sampled values of $\theta$. The two parameters were selected through a process of manual tuning, whereby various combinations were tested and assessed for their performance on a subset of problems within the dataset. The selection was guided by visual inspection of the resulting fits, with consideration given to balancing model complexity and accuracy. To further validate the result of the filtering we have checked that the filter maintains the positive sign of the data and we manually reviewed the cases where the difference between the minimum of the filtered and unfilterd (raw) data is larger than a certain threshold. In these cases we checked if it would be appropriate to increase the degree of the polynomial or change the size of the window. Indeed, the smoothing sometimes changes the position of the minima in sharp valleys. 

Computing time also depend on ``external''
parameters, such as the machine used for the computation. To overcome the variability due to the hardware and, at the same time, improve the training, we normalize the data between zero and one, independently for each subset of outputs in our dataset defined by fixing a matrix $\mathrm A$. In this way, we require only that the computations concerning the same $\mathrm A$ are carried on the same machine. Once the multigrid solver is used in parallel, another external factor that may influence the computational cost is the number of CPU cores. We can easily treat this variable as another parameter of the AMG method, however, as mentioned before, this extension is left for future works. The computations were carried out (in serial) on various nodes of the HPC cluster at MOX with processor Intel${}^\textnormal{\tiny\textregistered}$ Xeon${}^\textnormal{\tiny\textregistered}$ Gold 6238R.

\subsection{Applying the pooling operation to large sparse matrices}\label{sec:pooling}
The first step of our algorithm is to transform the large sparse matrix $\mathrm A$ into a small tensor (that can be seen as a multi-channel image) so that can be cheaply processed by the ANN. Drawing on the well-established benefits of pooling, which include translation invariance, computational efficiency, and hierarchical representation (see \ref{apx:dl} for details). We aim at extending its application to the sparse matrix A. Indeed, $\mathrm A$ can be seen as a large gray scale image and we want to apply the pooling operation (before the training) to reduce the computational cost associated with handling such a large amount of data, and to prune unimportant features. Our hypothesis is that not all the information contained in $\mathrm A$ is needed for the regression problem at hand. We call $\mathrm V = \texttt{pooling} (\mathrm A, m) \in \mathbb R^{m\times m \times 4}$ the result of the pooling of $\mathrm A$, where $m$ is a hyperparameter of the ANN that defines the size of $\mathrm V$. We report the details in Algorithm~\ref{a:pooling}. Here, we are assuming that the sparse matrix $\mathrm A$ is stored in coordinate list format, however the algorithm can be easily generalized to other sparse storage formats. The only difference with the pooling operation usually applied in CNNs is that instead of extracting the maximum, we extract the following features (always in a rectangular neighborhood): maximum of the positive part, maximum of the negative part, the sum, the number of non-zeros (\textit{nnz}) entries. The insight behind extracting these features is the following: the positive part, the negative part and the sum are relevant in the definition of the weights of the interpolation operator Eq.~(\ref{eq:interp_w_def}), the \textit{nnz} count is an indicator of the sparsity of the neighborhood. Adding other features could be helpful, but we found this to be a good compromise between accuracy and computational cost. Moreover, these feature preserve the sparsity of the matrix, that is, if we consider a rectangular neighborhood of $\mathrm A$ of only zeros the result will be zero.
The algorithm has low complexity, namely $O(nnz)$, that in the case of FE means $O(p\cdot n_1)$. Moreover, we have supported experimentally that its elapsed CPU time is negligible with respect to the one to solve the linear system. This is a necessary condition for this algorithm to be worthwhile. We remark that Algorithm~\ref{a:pooling} could be easily generalized to work in parallel.

\begin{algorithm}[t]
    \caption{Pooling algorithm $\mathrm V=\texttt{pooling}(\mathrm{A},m)$}
    \label{a:pooling}
    \begin{algorithmic}[1]
        \STATE access $\mathrm{A}$ in COO form and extract its: size $n_1$, \texttt{val}, \texttt{row}, \texttt{col}
        \STATE initialize $\mathrm{V}$ to an $m \times m \times 4$ dense tensor with all zero entries $v_{ijl} = 0$
        \STATE $q$ $\leftarrow$ $n_1 / m$, $p$ $\leftarrow$ $n_1$ \textnormal{mod} $m$, $t$ $\leftarrow$ $(q+1)p$
        \FOR{$k = 0$ \textbf{to} $\textnormal{\texttt{val.size()}}-1$}
            \STATE $i\,$ $\leftarrow$ \texttt{row}[$\,k\,$]$/(q+1)$ \textnormal{\textbf{if}} (\texttt{row}[$\,k\,$] $< t$) \textnormal{\textbf{else}} $($\texttt{row}[$\,k\,$] $- \,t)/q + p$
            \STATE $j$ $\leftarrow$ \texttt{col}[$\,k\,$]$/(q+1)$ \textnormal{\textbf{if}} (\texttt{col}[$\,k\,$] $< t$) \textnormal{\textbf{else}} $($\texttt{col}[$\,k\,$] $- \,t)/q + p$
            \STATE $v_{ij1}$ $\leftarrow$ {$\max\{\max\{0, \texttt{val}[\,k\,]\}, v_{ij1}\}$}
            \STATE $v_{ij2}$ $\leftarrow$ {$\max\{\max\{0, -\texttt{val}[\,k\,]\}, v_{ij2}\}$}
            \STATE $v_{ij3}$ $\leftarrow$ {$v_{ij3} + \texttt{val}[\,k\,]$}
            \STATE $v_{ij4}$ $\leftarrow$ {$v_{ij4} + 1$}            
        \ENDFOR
    \end{algorithmic}
\end{algorithm}

\subsection{Normalization}\label{sec:normalization}
Still $\mathrm V$, defined as before, cannot be directly used as input of an ANN. Normalization of data is a necessary step to assure fast convergence in neural networks. To this end we employ the following logarithmic normalization for each channel of the input $\mathrm V$ since it has been shown in \cite{antonietti2023accelerating} to outperform a classical linear normalization in this kind of applications:
\begin{equation}
    \hat v_{ij} = \frac{\log(|v_{ij}| + 1)}{\max_{i, j}|\log(|v_{ij}| + 1)|} \frac{v_{ij}}{|v_{ij}|}.
\label{eq:normalization}
\end{equation}
One of the main properties of this normalization is to map zero in zero, one could see this property as preserving the sparsity pattern of $\mathrm A$. In the Appendix~\ref{apx:pooling_visual} we show some examples of the combination of pooling and normalization for some matrices $\mathrm A$ and confront them with the relative sparsity pattern.

\subsection{The neural network architecture}\label{sec:ann-architecture}
Once we have $\hat{\mathrm{V}}$, we feed it in the neural network shown in Figure~\ref{fig:architecture}, which we are now going to describe. We employ an ANN $\mathscr F$ to predict, as target the computational cost $t$ and the square of the error committed by the ANN itself $\sigma^2$. The insight behind $\sigma$ is that, since the elapsed CPU time $t$ is polluted with noise due to the measurements, in this way we are able to know when the ANN is confident on its own prediction. 

The architecture of $\mathscr F$ is comprised by two components. The first one leverages a CNN to extract the relevant features from the matrix $\mathrm V$. These features are flattened in a intermediate dense layer, where they are concatenated with the other features $p, \log_2(n_1)$ and $\theta$ and fed into a dense feed forward network, which predicts the computational cost $t$. On the output layer we clip the prediction of the normalized computational cost between zero and one and use a softplus activation function for the variance estimate.

We employ as inputs of the ANN the vector formed by the pooling of the matrix a $\hat{\mathrm{V}}$, the normalized result of $\texttt{pooling}(\mathrm A, m)$, where the size $m$ is a
hyperparameter that we tune (see Section~\ref{sec:hyperparameters-study}), the FEM degree $p$, the logarithm of the size of A $\log_2(n_1)$ and the strong threshold $\theta$. Namely, 
\begin{linenomath}\begin{equation*}
\mathscr{F}(\mathrm V, p, \log_2(n_1), \theta; \boldsymbol{\gamma}) = (\tilde{t}, \tilde{\sigma})
\end{equation*}\end{linenomath}
and we optimize (the superscript $i$ indicates the $i$-th sample of the dataset of size $T$):
\begin{equation}\label{eq:loss-with-err}
\min_{\boldsymbol{\gamma}} \frac{1}{T}\sum_{i=1}^T \textnormal{MSE}(t^i, \tilde{t}^i) + \textnormal{MSE}(\tilde{\sigma^i}^2, (t^i - \tilde{t}^i)^2).
\end{equation}
Eq.~(\ref{eq:loss-with-err}) could benefit from introducing weights for the individual loss terms, however, the order of magnitude of the two terms is similar, so we prefer to avoid adding an additional hyperparameter.
Let us remark that $p$ is not needed as input of the neural network. Indeed, the information about $p$ is embedded into the matrix $\mathrm A$: empirical results show that if we add $p$ as input of the network, we need a lower number of epochs to reach the same loss, but it is still possible to effectively train the network without $p$. Hence, our algorithm is not limited to problems stemming from FE discretizations.

Given the matrix $\mathrm A$ associated to the linear system we want to solve, then our algorithm prescribes the default literature value of the strong threshold $\theta$ every time the weighted average variance $\hat{\sigma}$ of the map $\mathrm A \mapsto \theta$,
\begin{equation}\label{eq:sigma-hat}
\hat{\sigma} = \frac{1}{181}\sum_{j=10}^{190}(1 - \tilde{t}^j)\tilde{\sigma}^j, \textnormal{  where  } (\tilde{t}^j, \tilde{\sigma}^j) = \mathscr{F}(\mathrm V, p, \log_2(n_1), \frac{j}{200}; \boldsymbol{\gamma})
\end{equation}
is greater than a certain threshold $\bar \sigma$. In Eq.~(\ref{eq:sigma-hat}) we consider $\theta=\frac{j}{200} \in [0.05, 0.95]$ because we found that smaller or larger values of $\theta$ do not influence in a significative manner the performance of the AMG. This weighted average ensures that the variance of the prediction is more relevant if closer to the (expected) minima. The value $\bar \sigma$ is calibrated offline on the validation dataset once after the training of the ANN. For each matrix $\mathrm A^i$ in the dataset we compute the relative error indicator  $\hat{\sigma}^i$. Then we propose to choose $\bar{\sigma}$ as the ordinate of the elbow point of the sorted errors indicators $\hat{\sigma}^i$. On the other hand, if $\hat \sigma < \bar \sigma$, the algorithm prescribes as $\theta$ the value $\theta^*$ found by solving the optimization problem
\begin{equation}
   \theta^* = \argmin_{\, \theta \in (0,1]} \, \, (\mathscr F(\mathrm V, p, \log_2(n_1), \theta; \boldsymbol{\gamma}))_1.
   \label{eq:A-theta-map}
\end{equation}
However, we empirically found that by solving the discretization of the above with two hundred linearly spaced points, we have an estimate of $\theta^*$ that, for our aim, is indistinguishable from the solution of the continuous problem with gradient descent. There are three reasons to use this approach instead of predicting directly $\theta^*$:
\begin{enumerate}
    \item We can quantify the expected improvement with respect to using the standard literature value of $\theta$.
    \item We do not have to solve a new optimization problem just to add one sample to the dataset.
    \item From our experiments, this approach seems to be more robust with respect to directly predicting $\theta^*$.
\end{enumerate}

\subsection{The AMG-ANN algorithm}\label{sec:amg-ann-algo}
We now have all the ingredients to describe the proposed algorithm in detail. \paragraph{Offline phase} First, collect, smooth and normalize the data following the procedure outlined in Section \ref{sec:data-gather}. Then train the neural network with the collected samples. More details on how this is practically done are given in Section~\ref{sec:numerical-results}.

\paragraph{Online phase} It is detailed in the following, cf. also Algorithm~\ref{a:ann_amg}: as prerequisites the algorithm requires the matrix $\mathrm A$ of the linear system to be solved and the polynomial degree $p$ of the discretization.
\begin{itemize}
    \item[1] Compress $\mathrm A$ by applying the pooling to obtain $\mathrm V$ (Section~\ref{sec:pooling}).
    \item[2] Normalize $\mathrm V$ to obtain $\hat{\mathrm{V}}$ (Section~\ref{sec:normalization})
    \item[3-5] Get $\theta^*$ by finding the $\theta$ that fed in the ANN, together with $\hat{\mathrm{V}}$, $p$, $\log_2(n_1)$, minimizes the predicted computational cost (Section~\ref{sec:ann-architecture}).
    \item[6] Use $\theta^*$ as parameter of the AMG method.
\end{itemize}

\begin{algorithm}[t]
\caption{ANN-enhanced AMG \newline {$\mathbf{u} = \texttt{ANN\_AMG}(\mathbf{u}, \mathrm{A}, \mathbf{f}, \{(\mathrm{B}_1^{(j)}, \mathrm{B}_2^{(j)})\}_{j=k}^M, \nu_1, \nu_2, N_{max}, tol)$}}
\label{a:ann_amg}
\begin{algorithmic}[1]
    \STATE $\mathrm{V} \leftarrow  \texttt{pooling}(\mathrm{A}, m)$ (Algorithm~\ref{a:pooling})
    \STATE  normalize $\mathrm{V}$ by means of Eq.~(\ref{eq:normalization}), obtaining $\hat{\mathrm{V}}$

    \STATE  $(\tilde{t}^j, \tilde{\sigma}^j) \leftarrow \mathscr{F}(\hat{\mathrm{ V}}, p, \log_2(n_1), \frac{j}{200}; \boldsymbol{\gamma}), \quad j=10,11,...,190$
    \STATE $\hat{\sigma} \leftarrow \frac{1}{181}\sum_{j=10}^{190}(1 - \tilde{t}^j)\tilde{\sigma}^j$
    \STATE $\theta^* \leftarrow \argmin_{\tilde{t}^j} \frac{j}{200} \textnormal{\textbf{ if }} \hat{\sigma} > \bar{\sigma} \textnormal{\textbf{ else }} 0.5$
    \STATE  $\mathbf{u}  \leftarrow \texttt{AMG}(\mathbf{u}, \mathrm{A}, \mathbf{f}, \theta^*, \{(\mathrm{B}_1^{(j)}, \mathrm{B}_2^{(j)})\}_{j=k}^M, \nu_1, \nu_2, N_{max}, tol)$ (Algorithm~\ref{a:amg})
\end{algorithmic}    
\end{algorithm}
 
\subsubsection{Evaluating the performance of the algorithm}\label{sec:perf-indexes}
A small loss is not a good indicator of the performance of our algorithm. 
Indeed, the choice of $\theta$ is subordinate to the map $\mathrm A \mapsto \theta^*$ defined by Eq.~(\ref{eq:A-theta-map}). With this in mind, we introduce the following quantities of interest. Let $\mathrm A$ be fixed, and let:
\begin{itemize}
    \item $t_{\textnormal{ANN}}$ be the computational time of the AMG-ANN algorithm
    \item $t_{0.25}$ be the computational time of the AMG method for $\theta=0.25$
    \item $t_{\textnormal{MIN}}$ be the computational time of the AMG method with 
\begin{linenomath}\begin{equation*}
 \theta^*=\argmin_{\theta \in \textnormal{dataset for A}} t(\theta; \mathrm A).
\end{equation*}\end{linenomath}
\item $P = 1 - \frac{t_\textnormal{ANN}}{ t_{0.25}}$ be the performance index of the AMG-ANN algorithm
\item $P_\textnormal{MAX} = 1 - \frac{t_{MIN}}{ t_{0.25}}$ be the best performance of the AMG-ANN algorithm (according to the dataset).
\end{itemize}
Moreover, we can compound the quantities over different $\mathrm A$ and define $PB$ as the percentage of cases where $P \geq 0$, $P_m$ as the average of $P$ and $P_M$ as the median of $P$. The ratio $P/P_{MAX}$ gives a measure of how well the ANN has learned the data.

\section{Numerical Results}\label{sec:numerical-results}
One of the advantages of our algorithm is that it integrates seamlessly into pre-existing code. We carry out several numerical experiments using using deal.II \cite{dealii2019design, arndt2020deal} for the construction of the linear system and we rely on BoomerAMG of the library HYPRE \cite{falgout2002hypre} as the AMG solver. No changes are made to BoomerAMG, we use it as a black-box solver and only change the value of the strong threshold parameter $\theta$, all the other parameters (choice of pre- and post- smoother $\mathrm B$, number of pre- and post- smoothing cycles $\nu$, etc.) are left as default. We set the tolerance on the (relative) residual to $10^{-8}$.

Concerning the training of the ANN, we employ a 20-20-60 split of the dataset into training-validation-test, respectively. The splitting is done among problems, this means that a certain matrix $\mathrm A$ appears just in one of the three datasets. This entails that when we evaluate the algorithm on the test (or validation) dataset, it is making predictions on problems it has never seen before. Notice also that the test dataset is much larger than the other two: this makes sense only when the dataset is large enough so that the neural network can generalize well. Training has been performed using the Adam optimizer with default learning rate $10^{-3}$ and minibatch size equal to 32. Moreover, we found that applying a suitable learning rate schedule is key to accelerate the optimization. Namely, we employ a learning rate schedule that halves the learning rate with patience 15 epochs.

The convolutional part of the network is composed by one to three plain convolutional blocks each one ending with a max-pooling layer. For each test case, we tune the number of convolutional layers, the size of the filter, the number of filters the number of dense layers and the width of the dense layers.

\subsection{Test Case 1: Highly heterogeneous diffusion problem, unstructured Grids}\label{sec:tc2}
Let us consider the following parametric elliptic problem
\begin{equation}
\begin{cases}
-\text{div}(\mu\nabla u) = f, \quad &\text{in} \: \Omega,\\
u = 0, &\text{on} \: \partial \Omega,
\end{cases}
\label{eq:poisson}
\end{equation}
where $\mu \in L^\infty(\Omega)$ is a highly heterogeneous piecewise positive constant and $\Omega$ varies among four geometries: a simplex, a plate with a hole, a ball and a cylinder. We consider in total eight different discretizations of the four domains and their nested refinements. We show the coarsest meshes in the first column of Table~\ref{tab:meshes}. The diffusion coefficient $\mu$ is a highly heterogeneous piecewise constant that is conforming with the mesh. Namely, $\mu$ has a different value equal to $10^{\varepsilon_i}$ on the $i$-th cell of one of the coarsening of the mesh, where $\varepsilon_i$ is randomly chosen between $[0, \varepsilon_{MAX}]$ and $\varepsilon_{MAX}=1, 3, 10$. The problem is discretized by means of continuous FE of order $p=1,2,3$. In order to test the stability of the algorithm we also use four different DoFs (degrees of freedom) numbering. When applying the renumbering, the underlying sparsity graph of the matrix $\mathrm A$ and thus the optimal value of $\theta$ stays the same, but the pooling $\mathrm V$ changes, so the ANN should learn to be invariant with respect to these changes. Namely we employed: the default dealii numbering, the Cuthill-McKee renumbering method implemented by dealii and by the Boost library \cite{siek2001boost} and King renumbering. We solve each problem with 37 values of $\theta$ with equispaced values between 0.05 and 0.95. In this way, we build a dataset counting 5873 different configurations and thus containing 217301 samples.

\begin{table}[t]
\small
  \centering
  \begin{tabular}{ | c | m{2.1cm} | m{7.5cm} | }
    \hline
    Picture & Geometry & Description \\ \hline
    \begin{minipage}{.1\textwidth}
      \includegraphics[width=\linewidth]{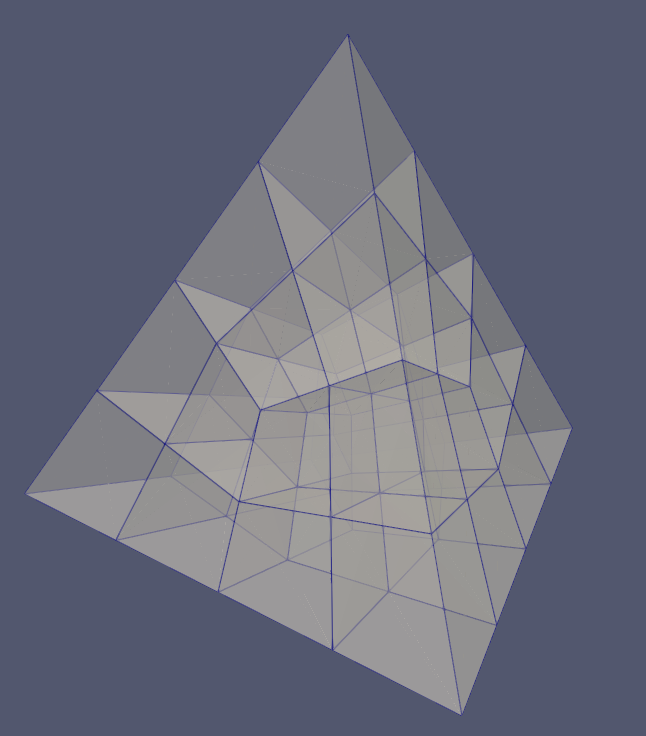}
    \end{minipage}
    &
    Simplex
    & 
    Tetrahedron inscribed in the unit ball centered in $(0, 0, 0)$.
    \\ \hline
    \begin{minipage}{.1\textwidth}
      \includegraphics[width=\linewidth]{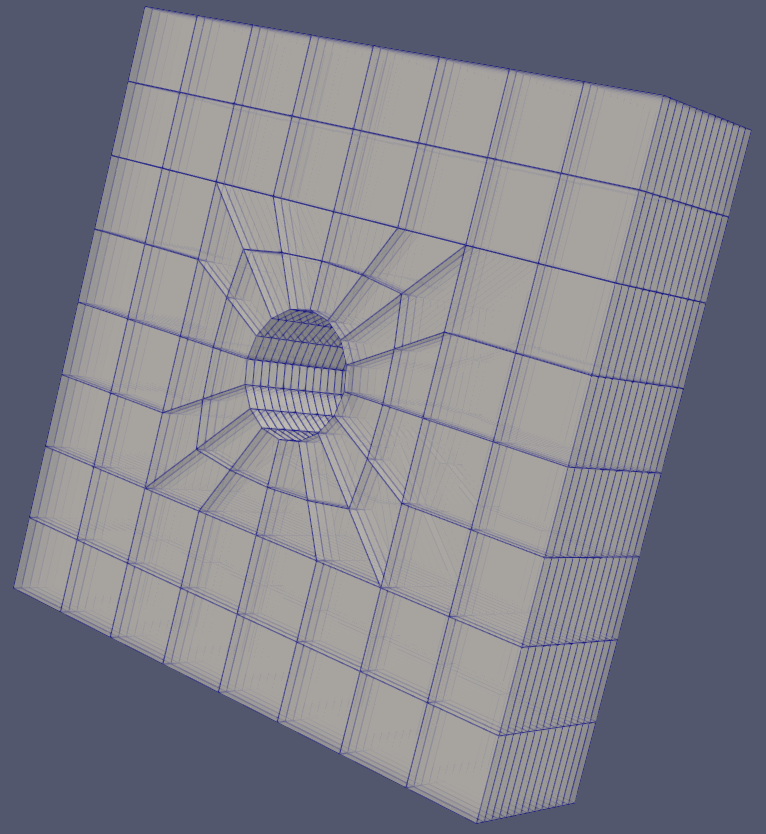}
    \end{minipage}
    &
    Plate with hole
    & 
    The $(-2, 2)\times(-2, 2)\times(-\frac{1}{2}, \frac{1}{2})$ plate with a centered hole of radius 0.4 along the z axis. Discretizations with two or eight slices in the z-direction are considered.
    \\ \hline
    \begin{minipage}{.1\textwidth}
      \includegraphics[width=\linewidth]{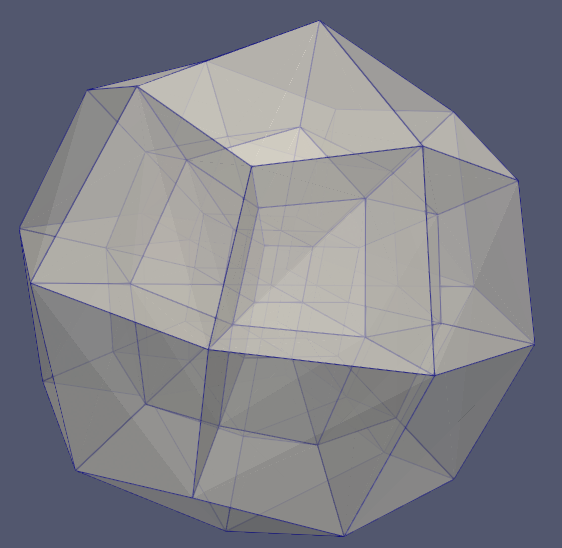}
    \end{minipage}
    &
    Ball
    & 
    Unit ball centered in zero.
    \\ \hline
    \begin{minipage}{.1\textwidth}
      \includegraphics[width=\linewidth]{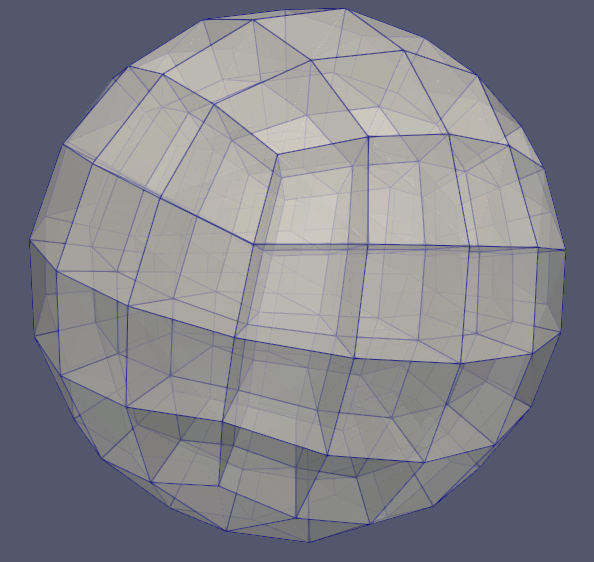}
    \end{minipage}
    &
    Balanced ball
    & 
    A variant of the unit ball centered in zero that has a better balance between the size of the cells around the outer curved boundaries and the cell in the interior.
    \\ \hline
    \begin{minipage}{.1\textwidth}
      \includegraphics[width=\linewidth]{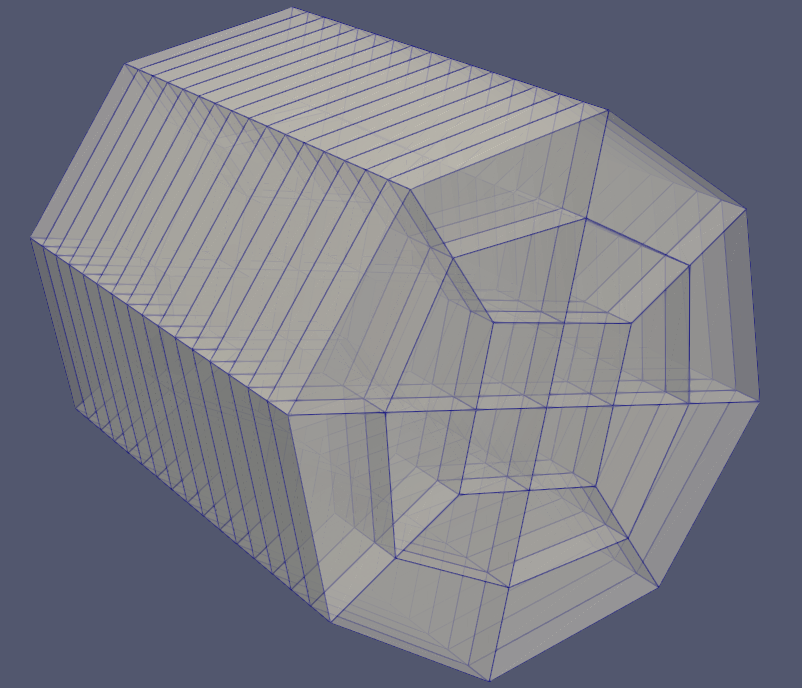}
    \end{minipage}
    &
    Cylinder
    & 
    Cylinder of unit radius and height two centered in zero. We consider three different discretizations with one, two or eight slices along its height. 
    \\ \hline
    \begin{minipage}{.1\textwidth}
      \includegraphics[width=\linewidth]{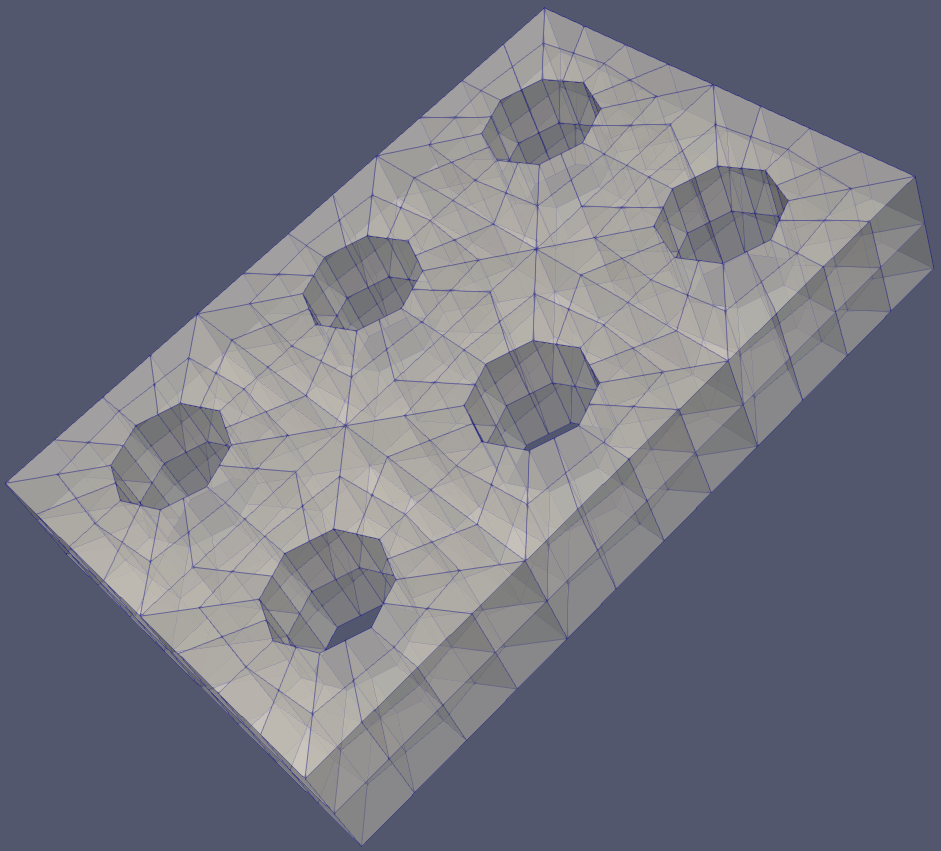}
    \end{minipage}
    &
    Holes
    & 
    $3\times2$ replica of the ``Plate with hole'' geometry on the $x$ and $y$. Only test dataset. 
    \\ \hline
    \begin{minipage}{.1\textwidth}
      \includegraphics[width=\linewidth]{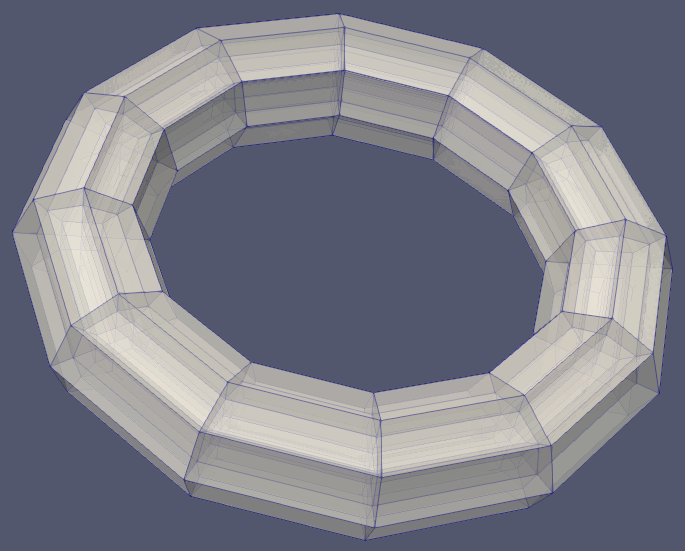}
    \end{minipage}
    &
    Torus
    & 
     Torus of circle radius 2 and inner radius $\frac{1}{2}$. The axis of the torus is the $y$-axis. Only test dataset.
    \\ \hline
  \end{tabular}
  \caption{Test case 1: list of meshes used in the dataset.}\label{tab:meshes}
\end{table}

\subsubsection{Study of hyperparameters and pooling features}\label{sec:hyperparameters-study}
\begin{figure}[t]
    \makebox[\textwidth][c]{\includegraphics[trim={0cm 0.5cm 0cm 0cm},clip,width=0.9\textwidth]{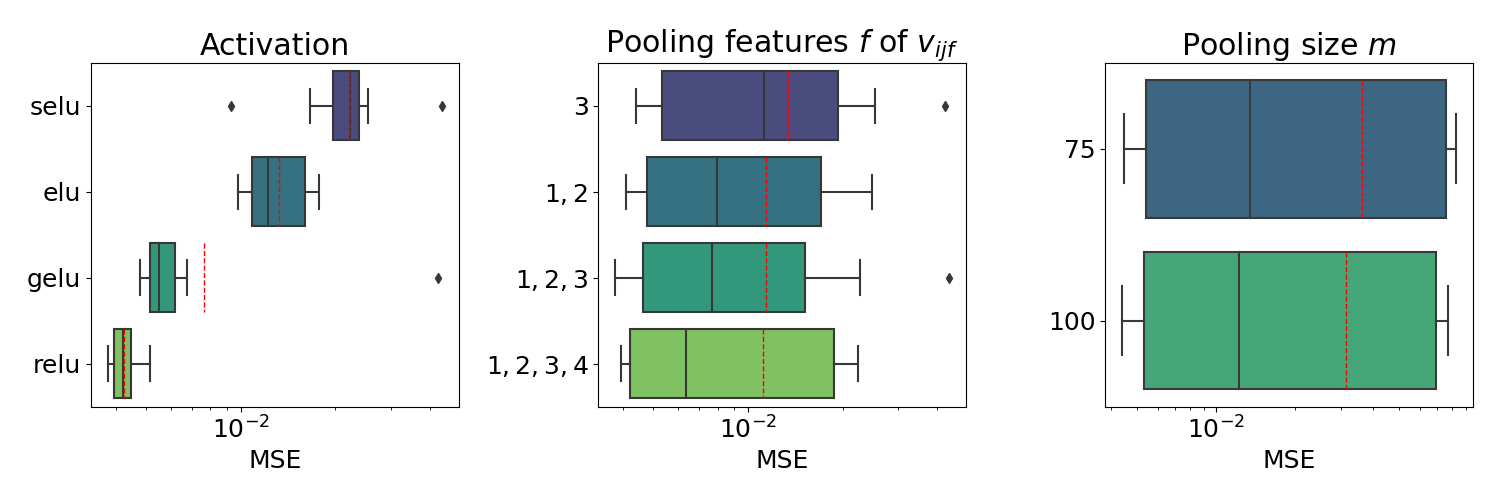}}
    \caption{Boxplots of the impact of the choice of the activation function (\textit{left}), the pooling features (\textit{center}) and the pooling size $m$ (\textit{right}), on the validation loss (MSE) of an ANN trained on a subsample of the dataset described in Section~\ref{sec:tc2}. The red line represents the mean. }
    \label{fig:hyper-params-study}
\end{figure}

In this section we analyze the influence of some hyperparameters on the final loss when we employ a subset of the dataset. This is a preliminary analysis that we made for each test case in order to choose the best value of $m$ (the size of the pooling $\mathrm V$), which of the features of $\mathrm V$ are relevant (that is which of the four layers $f$ of $v_{ijf}$ should we employ as input of the CNN) and the activation function. Namely we performed a grid search over different network architectures. Figure~\ref{fig:hyper-params-study} shows the results. The ReLU activation function consistently reaches the smallest loss, hence we employ ReLU activation functions in conjunction with He initialization \cite{he2016deep}. The ANN that uses all four the features of the pooling has the smallest loss, thus the whole tensor $\mathrm V$ is relevant for this task. The pooling size $m=75$ gives the best trade-off between accuracy and computational cost.

\subsubsection{Calibration of the error threshold}
\begin{figure}[t]
        \centering
        \includegraphics[width=0.45\linewidth]{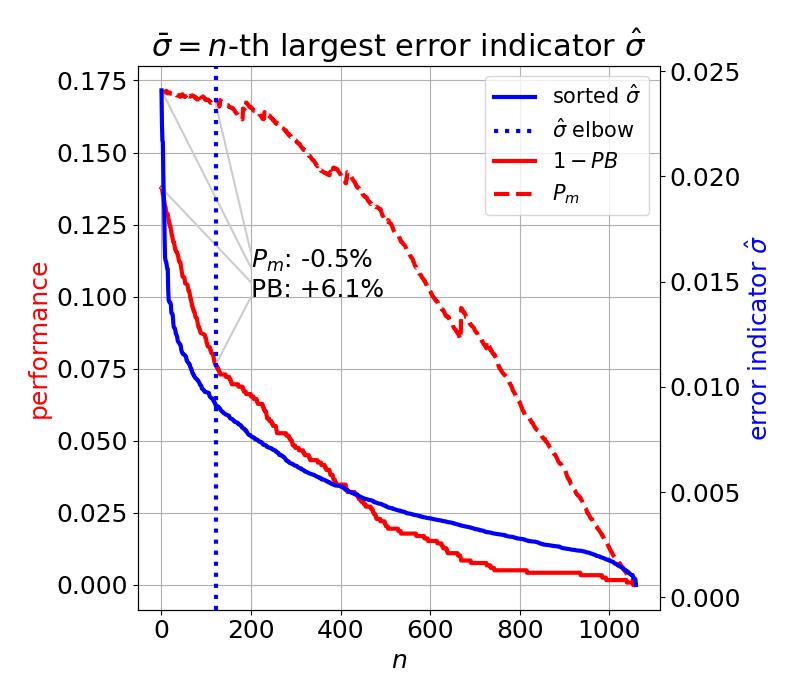}
        \caption{Test case 1: performance accuracy trade-off of the AMG-ANN algorithm depending on the choice of $\bar \sigma$. The algorithm is evaluated on the validation set.}
        \label{fig:sigma_bar_choice}
\end{figure}

\begin{figure}[t]
    \centering
    \includegraphics[trim={0.8cm 0cm 0.0cm 0cm},clip,width=0.8\linewidth]{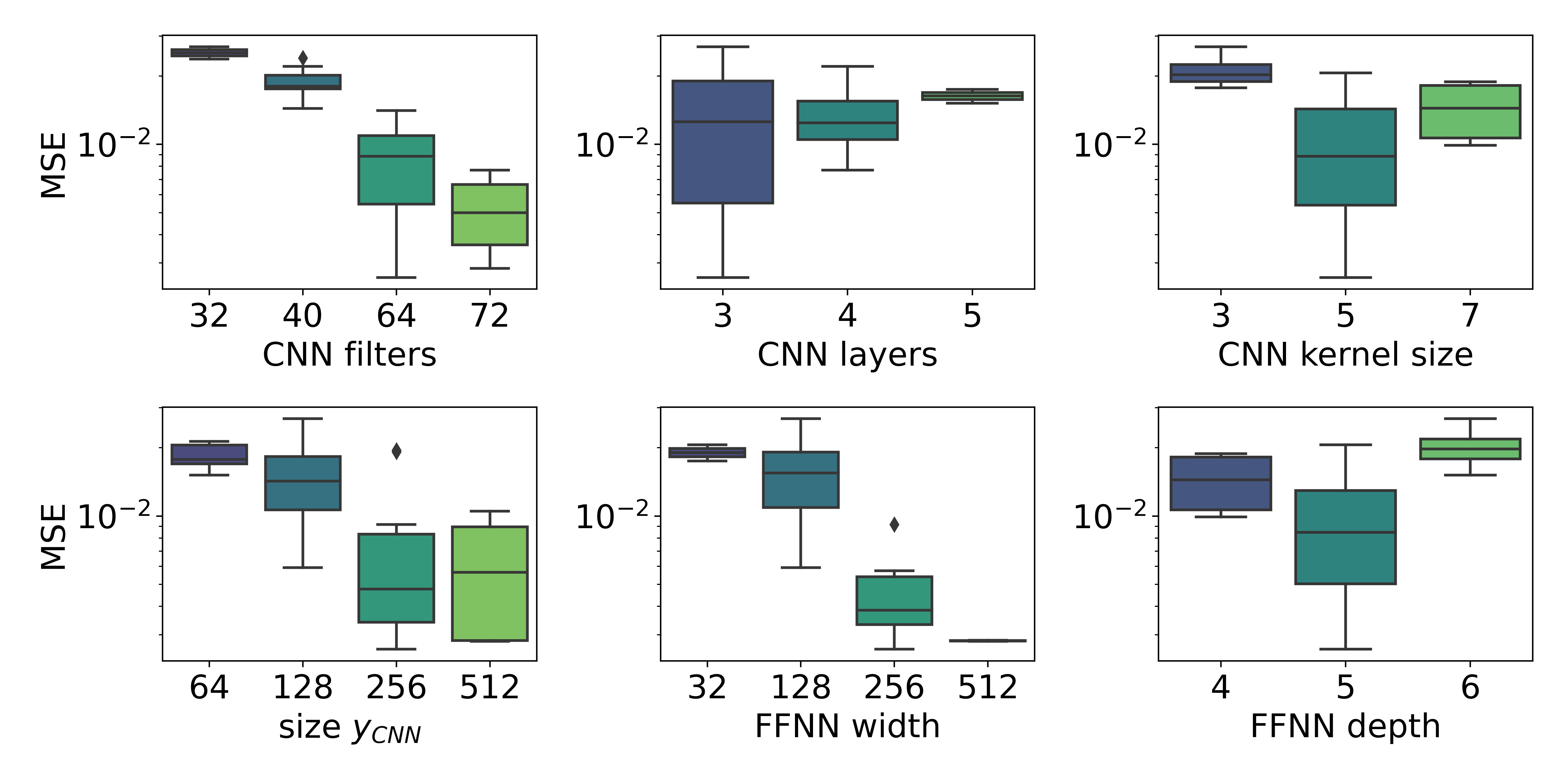}
    \caption{Test case 1: Box plot of the validation loss values after 80 epochs depending on the choice of the following hyperparameters: the number of filters in the CNN layers, the number of CNN layers in each CNN block, the kernel size in the CNN layers, the size of the output of the CNN, the width and depth of the dense FFNN.}
    \label{fig:hyper-fine-tune}
\end{figure}

\begin{figure}[t]
    \centering
        \centering
        \includegraphics[trim={0cm 0cm 0cm 0cm},clip,width=0.4\linewidth]{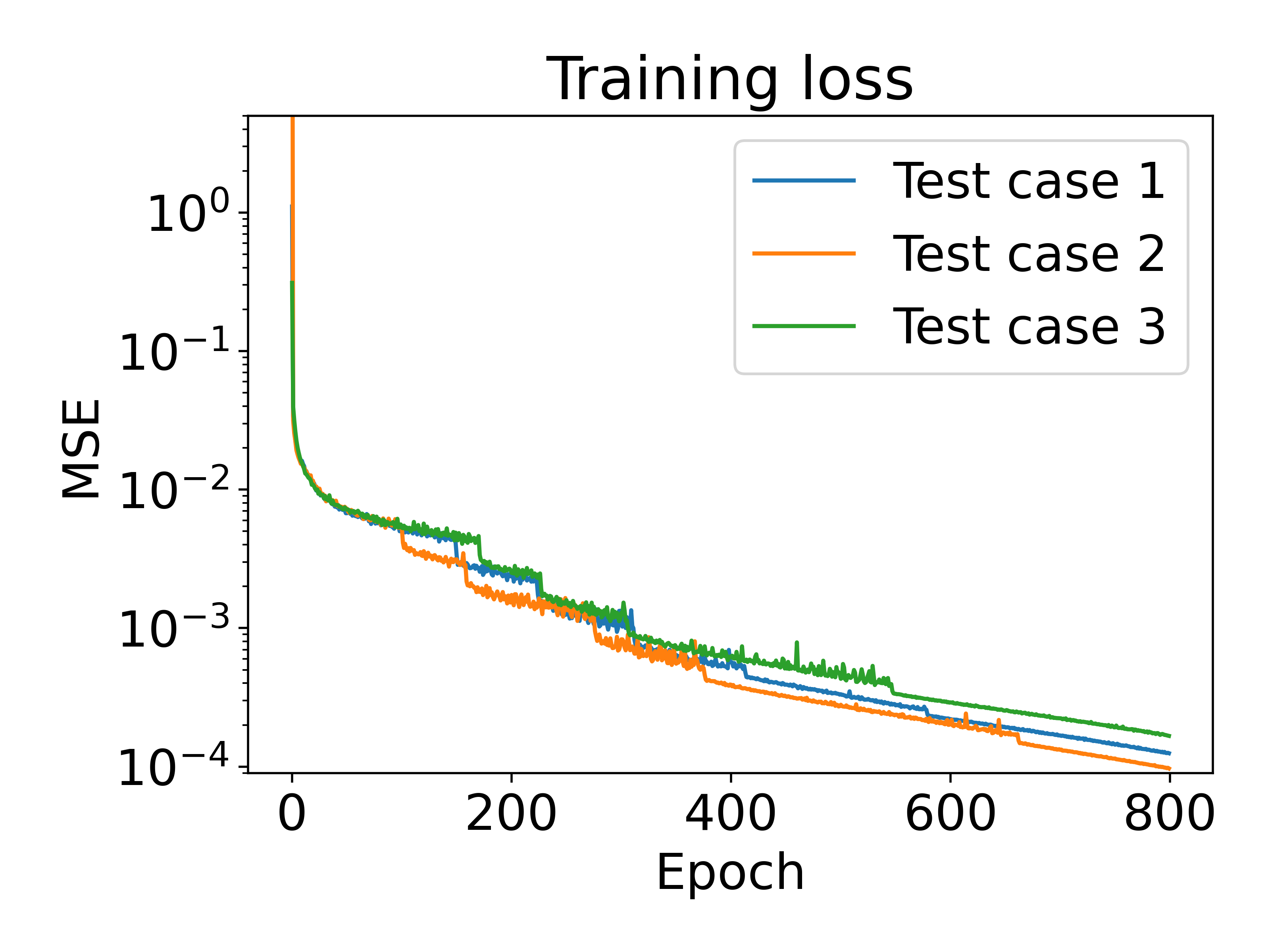}
        \caption{History of the training loss of the fine-tuned ANN.}
        \label{fig:loss-history}
   
\end{figure}

\begin{figure}[t]
\setlength{\tabcolsep}{0pt}
    \centering
        \begin{tabular}{c}
             \includegraphics[trim={0cm 0cm 0cm 0cm},clip,width=0.9\linewidth]{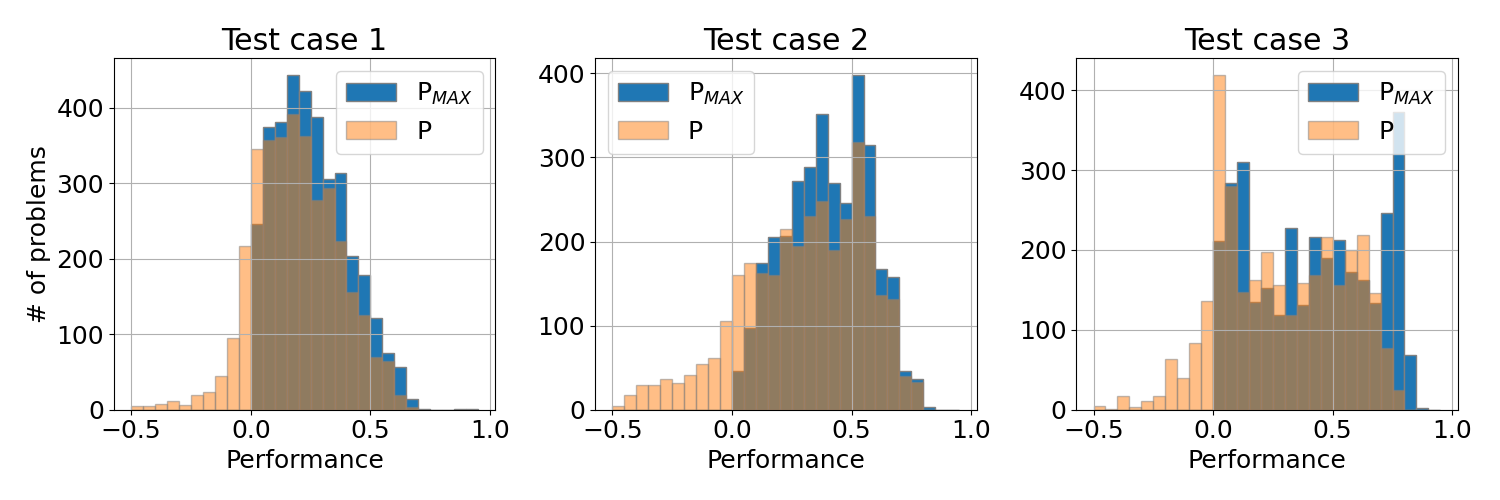}\\[-4pt]
             \includegraphics[trim={0cm 0.5cm 0cm 0cm},clip,width=0.9\linewidth]{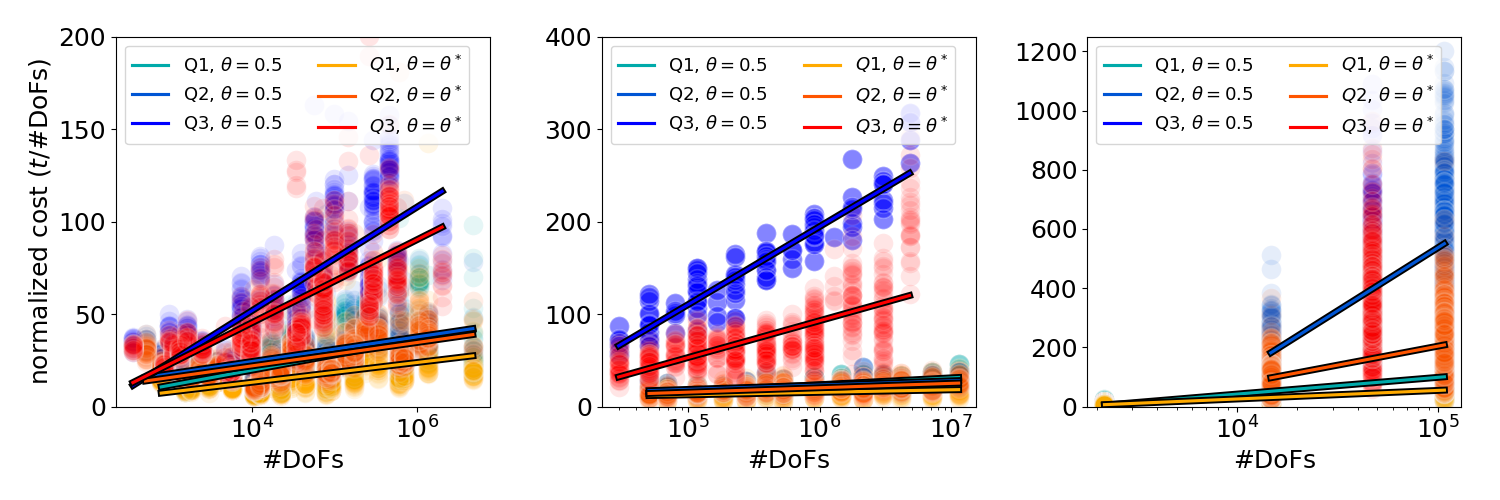}\\[-4pt]
        \end{tabular}%
    \caption{(\textit{Top}) Performance gain $P$ of the AMG-ANN algorithm and maximum theoretical performance $P_{MAX}$ for the three test cases. (\textit{Bottom}) Scatter plot with linear regression of the DoF normalized computational cost of the AMG-ANN algorithm with respect to the number of DoFs and the polynomial degree $p$. In tones of red the AMG-ANN algorithm and in tones of blue the standard AMG method.\label{fig:perf-and-cost}}
\end{figure}

\begin{figure}[t]
        \centering
        \includegraphics[trim={0cm 1cm 0cm 1cm},clip,width=0.45\linewidth]{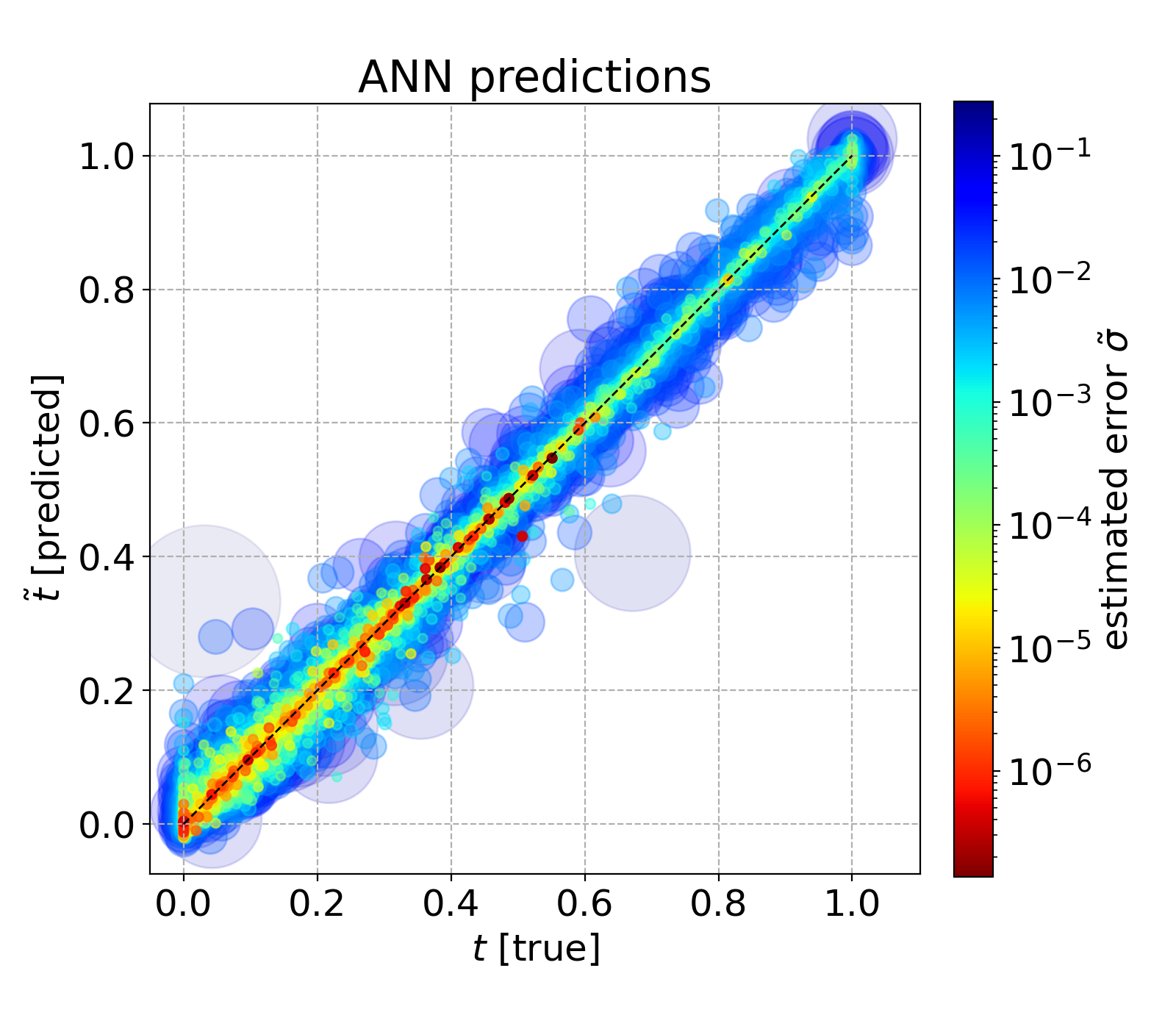}
        \caption{Test case 1: Prediction of the ANN. Color and size are proportional to the estimated error $\tilde{\sigma}$.}
        \label{fig:variance_prediction_scatter}
\end{figure}

We discuss how the choice of $\bar \sigma$ is made and how it impacts the performance of our algorithm. In Figure~\ref{fig:sigma_bar_choice} we plot the accuracy (PB) and average time reduction ($P_m$) when choosing as $\bar \sigma$ the $n$-th largest $\hat \sigma$ of the validation dataset. Notice there is an accuracy-performance trade-off, indeed, as $\bar \sigma$ tends to zero, $PB$ tends to one but at the same time $P_m$ tends to zero. Moreover, the fact that $P_m$ is not strictly decreasing means that the error committed by the algorithm is not just due to noisy measurements. However, by using as $\bar \sigma$ the elbow of the ordered $\hat \sigma$ we have a gain of 6.1\% in terms of accuracy meanwhile the mean performance decreases of only 0.5\%.

\begin{figure}[t]
    \includegraphics[trim={0 0cm 0 0cm},clip,width=0.95\textwidth]{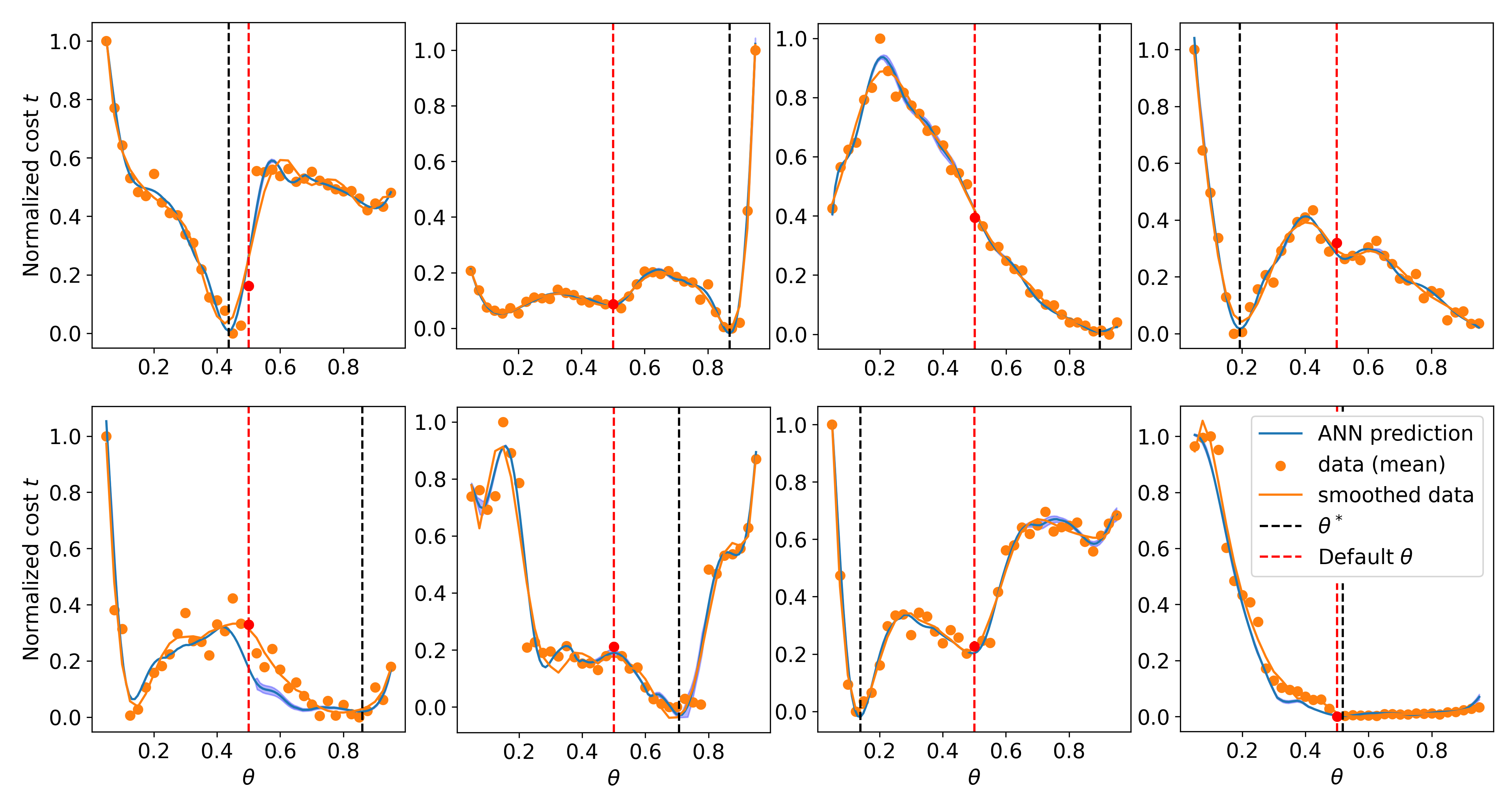}
    \caption{Test case 1: Predictions made by the ANN for some exemplary values of the diffusion coefficient and domain shape. On the y-axis the normalized computational cost $t$, on the x-axis the value of the strong threshold parameter $\theta$. The dots represent the raw measurements of the computational time.}
    \label{fig:predictions-examples}
\end{figure}

\begin{figure}[t]
    \makebox[\textwidth][c]{\includegraphics[trim={0cm 0cm 0.0cm 0cm},clip,width=0.712\linewidth]{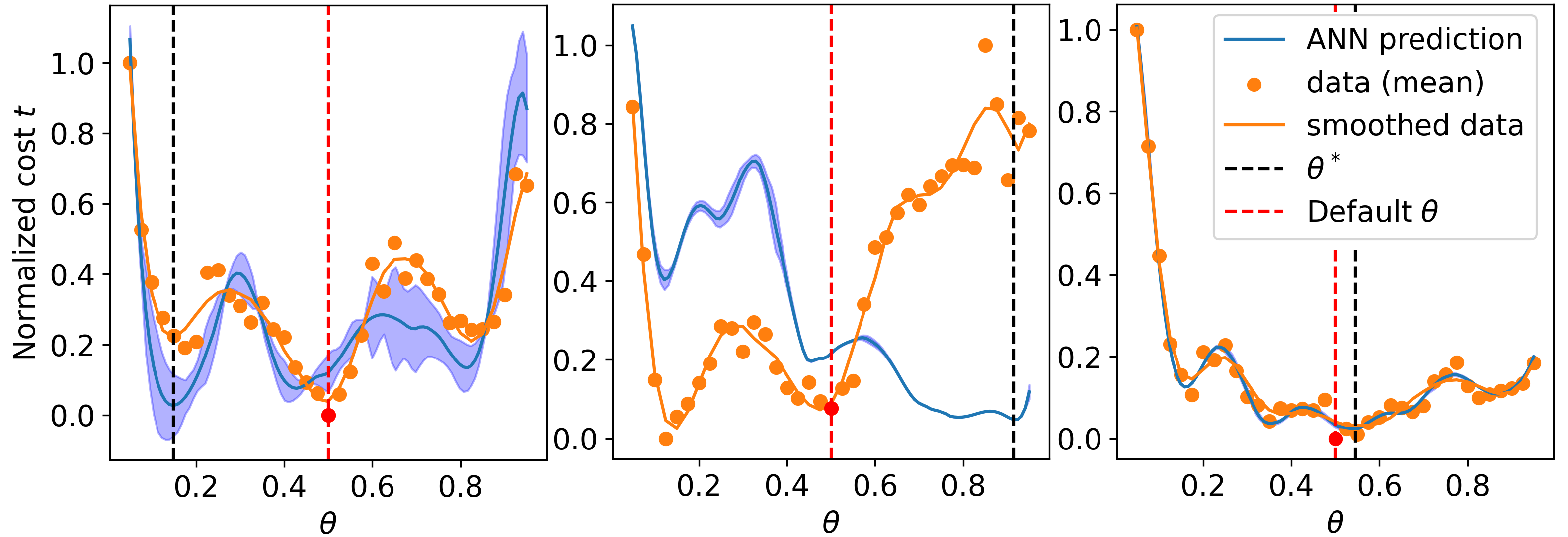}}
    \caption{Test case 1: three representative cases of predictions of the normalized computational cost $t$ made by the ANN that leads to sub-optimal values of $\theta^*$. The blue area is the ANN error estimate $\tilde\sigma$. }
    \label{fig:error-types}
\end{figure}

\subsubsection{Evaluation of the algorithm}
\begin{figure}[t]
    \makebox[\textwidth][c]{\includegraphics[trim={4cm 2cm 0 0cm},clip,width=0.9\textwidth]{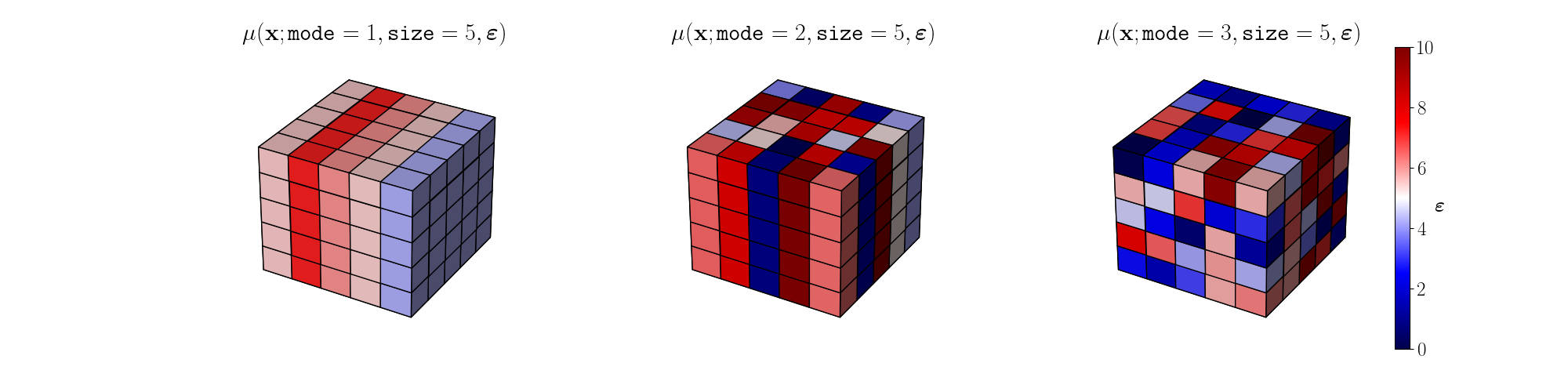}}
    \caption{Examples of the diffusion coefficient $\mu$ employed in Section~\ref{sec:tc1} for some choices of the vector $\boldsymbol \varepsilon$ depending on the \texttt{mode} parameter (one to three from left to right) and fixed \texttt{size}=5.}
    \label{fig:mu}
\end{figure}

\begin{figure}[t]
    \makebox[\textwidth][c]{\includegraphics[trim={0cm 0.5cm 0cm 0.5cm},clip,width=0.75\textwidth]{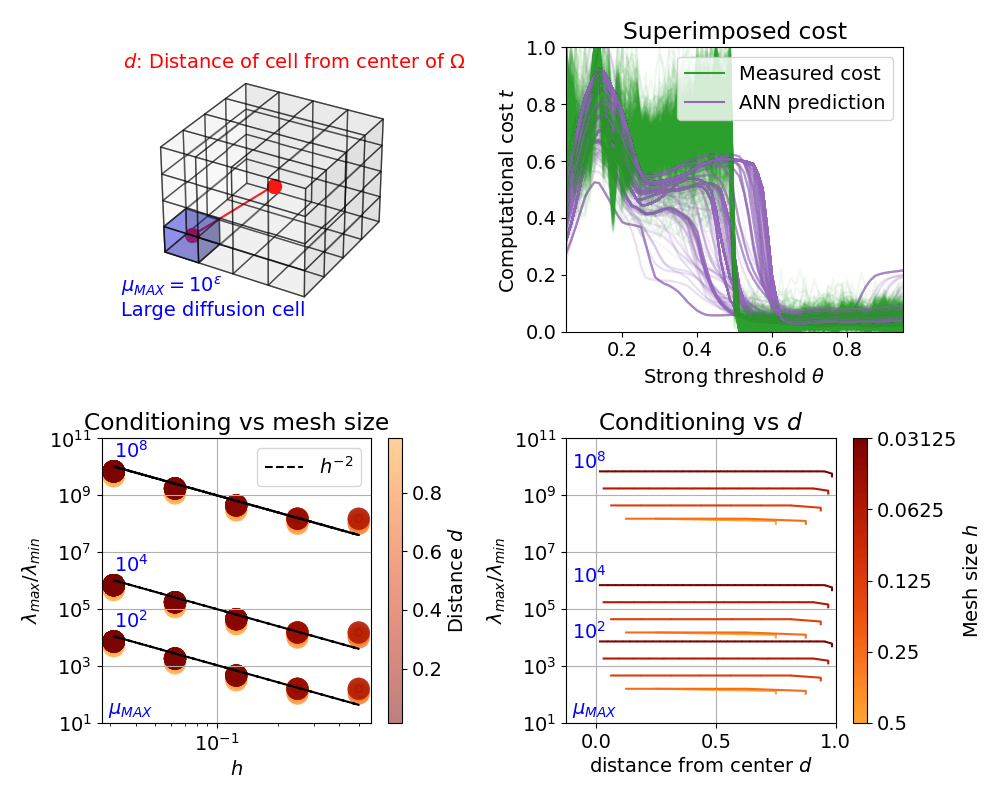}}
    \caption{Study of a special subset of problems of Test case 2. (\textit{Top-left}) We consider a uniform diffusion $\mu=1$ apart from one cell where $\mu=\mu_{MAX}=10^\varepsilon$ and $\varepsilon=2,4,8$. This cell has distance $d$ from the origin. (\textit{Top-right}) Superimposed normalized computational cost $t$ needed to solve the system with the AMG method depending on the strong threshold parameter $\theta$ and relative predictions made by the ANN.  (\textit{Bottom}) Conditioning number of the matrix $\mathrm A$ depending on the distance $d$ and the mesh size $h$.}
    \label{fig:spectrum}
\end{figure}

Starting from the architecture we found in Section~\ref{sec:hyperparameters-study}, we fine tuned the number of filters, the kernel size of the CNN and the width and depth of the dense part. Our strategy is to use a grid-search on a fine space with a limited number of epochs and then we increase the number of epochs while discarding the combinations of hyperparameters that have the largest validation loss. This process is repeated until we have one combination of hyperparameters. Figure~\ref{fig:hyper-fine-tune} shows an example of the loss depending on the hyperparameters after 80 epochs. We report in Table~\ref{tab:architecture} the final choice of hyperparameters. The ANN has a (test) loss of $2.89 \cdot 10^{-4}$, Figure~\ref{fig:loss-history} shows the history of the training loss. The number of iterations needed reach convergence is between 6 and 155 with $\theta=0.5$ and between 6 and 20 with $\theta=\theta^*$. The performance of our algorithm is evaluated by the indexes described in Section~\ref{sec:perf-indexes} on the test dataset (see Table~\ref{tab:perf_summary}) and by visualizing the scaling of the AMG with respect to the number of DoFs (see Figure~\ref{fig:perf-and-cost}). Moreover, we remark that the scaling shown in the bottom row of Figure~\ref{fig:perf-and-cost}, ideally, should be constant. Indeed, this is true when we apply the AMG to problems where $\max_\Omega \mu / \min_\Omega \mu \sim 1$, however in the cases that we consider is not possible to achieve this. In any case, the AMG-ANN is able to achieve a much better scaling with respect to the not tuned AMG version. In particular, the algorithm has an accuracy of 92.51\% and it produces in average a reduction of the computational cost of almost 16\%. Since in this case the median on $P_{MAX}$ is rather small (82\%), there is still margin of improvement for the reduction of the computational cost by tuning the architecture of the ANN.

\begin{table}[t]
\centering
\begin{tabular}{l|llll}
                       & Test case 1 & Test case 2 & Test case 3 &  \\
                       \hline
CNN filters            & 64          & 48          & 64          &  \\
CNN layers             & 3           & 3           & 3           &  \\
CNN kernel size        & 6           & 5           & 5           &  \\
$\mathbf y_{CNN}$ size & 256         & 256         & 256         &  \\
FNN width              & 512         & 256         & 512         &  \\
FNN depth              & 512         & 256         & 512         & 

\end{tabular}
\caption{Architecture of the ANN after the fine tune.}
    \label{tab:architecture}
\end{table}

\begin{table}[t]
    \centering
    \begin{tabular}{ ll | c cc c}
        Test case                             & $\bar \sigma$        & $PB$    & \multicolumn{2}{c}{$P$ (avg/median)} & $P/P_\textnormal{MAX}$ (median)\\ 
        
\hline
\multirow{2}{*}{1: Unstructured} & $\infty$             & 86.58\% & 16.73\% & 17.90\% & 85.54\%  \\
                                 & $\tilde\sigma$ elbow & 92.51\% & 15.91\% & 15.46\% & -  \\
\hline
\multirow{2}{*}{1: Holes} & $\infty$             & 80.82\% & 15.09\% & 18.02\% & 76.06\%   \\
                          & $\tilde\sigma$ elbow & 86.01\% & 15.13\% & 17.76\% & -   \\
\hline
\multirow{2}{*}{1: Torus} & $\infty$             & 63.88\% & 0.85\% & 5.69\% & 45.72\%   \\
                          & $\tilde\sigma$ elbow & 65.08\% & 0.01\% & 3.30\% & - \\
\hline
\multirow{2}{*}{2: Structured}   & $\infty$             & 87.97\% & 30.51\% & 32.14\% & 93.50\%  \\
                                 & $\tilde\sigma$ elbow & 93.03\% & 29.54\% & 19.42\% & -  \\
\hline
\multirow{2}{*}{3: Elasticity} & $\infty$             & 83.72\% & 24.01\% & 21.40\% & 88.59\%  \\
                               & $\tilde\sigma$ elbow & 90.71\% & 22.60\% & 18.54\% & -  \\
\hline
\multirow{2}{*}{3: Pretrain} & $\infty$             & 85.76\% & 23.99\% & 23.79\% & 88.49\%  \\
                               & $\tilde\sigma$ elbow & 92.39\% & 22.17\% & 22.31\% & -  \\
    \end{tabular}
    \caption{Evaluation of the performance of the AMG-ANN algorithm for each test case depending on the choice of the threshold error $\bar\sigma$.}
    \label{tab:perf_summary}
\end{table}

Figure~\ref{fig:variance_prediction_scatter} shows all the predictions made by the trained ANN, notice that the largest points are the furthest from the bisector, thus the ANN is successfully predicting where it is inaccurate.
In Figure~\ref{fig:predictions-examples} we highlight some relevant predictions. Notice that the map $\theta \mapsto t$ exhibits many different complex patterns and that $\theta^*$ varies depending on the problem considered. Hence, it is not possible to find a fixed value $\theta^*$ that works well for all the problems. In Figure~\ref{fig:error-types}, we showcase relevant cases of when the algorithm is not accurate, in particular, we categorize the errors of the algorithm into three classes (from left to right):

\begin{itemize}[topsep=0.1cm,itemsep=-0.5ex,partopsep=1ex,parsep=1ex]
    \item Noisy measurements: the error is due to the measurement of $t$. As mentioned before, we reduce this error by repeating the measurements and regularizing data. Moreover, $\tilde{\sigma}$ is often a good indicator of a large error in these cases.
    \item Generalization error: this takes place when the network is not able to generalize beyond the training set. You can notice that in this case the prediction is completely wrong.
    \item ANN error: in this category we classify the errors due to the  approximation capacity of the ANN, the optimization error and the loss of information occurring during the pooling of A. These type of errors can be mitigated by hyperparameters tuning. However, there are some limitations, for instance, even if we would like $m$ to be large, in order to lose the least amount of information possible, the cost of the evaluation of $\mathscr{F}$ depends quadratically on it, so it must stay small. Moreover, the error estimate $\tilde{\sigma}$ is affected by these kind of errors itself, and it may not be of use in these cases.
   
\end{itemize}

\subsubsection{Evaluation on unseen domains}\label{sec:unseen-domains}
We tested the ANN-AMG algorithm by applying it to problems with a diffusion coefficient $\mu$ that was not present in the training dataset, but on a known geometry. To prove the generalization capability of the algorithm, we now test it using domains and diffusion coefficients it has never seen before. The two domains we employ are represented in the two last rows of Table~\ref{tab:meshes}: a replication of the plate with a hole in the $x$ and $y$ coordinates and a torus. We build the dataset exactly as done before, in total it counts 445 different problems. The performance on the ``Holes'' geometry is similar to a known geometry: $PB=86\%$ and $P_M=18\%$. Thus, we have shown that the algorithm is able to generalize also on domains with a different topology w.r.t.\ the ones in the training dataset. This result is linked to the fact that a (small) subset of the considered geometry (the plate with the hole) is present in the training dataset. Indeed, the results on the Torus are much worse: $PB=64\%$ and $P_M=6\%$. The algorithm is still making predictions that are significantly better than any prediction made ``at random'', however this challenging mesh shows the limitation of our approach. We show the details about the performance in Table~\ref{tab:perf_summary}. 

\subsection{Test Case 2: Highly heterogeneous diffusion problem, structured grids}\label{sec:tc1}

\begin{algorithm}[t]
    \caption{Diffusion $\mu$ in a given point $\mathbf{x} \in \Omega = (-1, 1)^3$, having fixed the \texttt{mode}, \texttt{size} and $\boldsymbol{\varepsilon} \in \mathbb{R}^{\texttt{size}^{(\texttt{mode})}}$. \newline $\mu = \mu(\mathbf{x}; \texttt{mode}, \texttt{size}, \boldsymbol \varepsilon)$}
    \label{a:mu}
    \begin{algorithmic}[1]
        \STATE $j$ $\leftarrow$ $1$
        \FOR{$i = 0$ \textbf{to} \texttt{mode}}
            \STATE $j$ $\leftarrow$ $j + \lfloor ((\mathbf x)_i + 1) \texttt{size} / 2 \rfloor \, \texttt{size}^{i-1}$
        \ENDFOR
        \STATE $\mu$ $\leftarrow$ $10^{(\boldsymbol{\varepsilon})_j}$
    \end{algorithmic}
\end{algorithm}

Let us consider the same elliptic problem (\ref{eq:poisson}) of before. However, here we have $\Omega = (-1, 1)^3$ and $\mu \in L^\infty(\Omega)$ is a highly heterogeneous piecewise positive constant defined differently from before. Namely, $\mu$ is defined by means of Algorithm~\ref{a:mu}, where $\texttt{mode} = 1, 2, 3$ defines if the pattern is either made of slices, lines or is checkerboard like, $\texttt{size} = 1, 2, ...$ defines how many times the pattern repeats, and the vector $\boldsymbol \varepsilon \in \mathbb{R}^{\texttt{size}^{(\texttt{mode})}}$ defines the value of $\mu$. Figure~\ref{fig:mu} shows a representation of $\mu$ for $\texttt{size} = 5$ and  $\texttt{mode} = 1, 2, 3$. The problem is discretized by means of continuous FE of order $p$ on nested Cartesian meshes. The dataset is built by varying $p=1,2,3$, $\texttt{mode}=1,2,3$, $\texttt{size}=2,3,...,10$, the mesh size $h$ between $\frac{\sqrt[3]{2}}{\texttt{size}}$ and $\frac{\sqrt[3]{2}}{\texttt{size} 2^{8-p}}$, and choosing at random the components of the vector $\boldsymbol{\varepsilon}$ between 0 and $\varepsilon_{MAX}$, where $\varepsilon_{MAX}=1, 3, 10$. We employ again four different DoFs numbering to enhance the dataset. We solve each problem with 37 values of $\theta$ equispaced between 0.05 and 0.95. The dataset we built counts 5471 different problems and thus contains 202427 samples. 


We started by using an architecture with hyperparameters chosen as described in Section~\ref{sec:hyperparameters-study} and then tuned the number of filters and the size of the kernel of the CNN and the width and depth of the dense part of the ANN. 
We report in Table~\ref{tab:architecture} the final choice of hyperparameters. We obtain an ANN with a (test) loss of $1.55 \cdot 10^{-4}$, Figure~\ref{fig:loss-history} shows the history of the training loss. Details about its performance and its scaling are shown in Figure~\ref{fig:perf-and-cost} and Table~\ref{tab:perf_summary}. In particular, our algorithm has a 93\% accuracy with an average reduction of the CPU time of almost 30\%, and in half of the cases you can expect a reduction greater than 32\%. 
The number of iterations needed reach convergence is between 7 and 18 with $\theta=0.5$ and between 7 and 9 with $\theta=\theta^*$.

\subsubsection{Study of a well-known family of problems}\label{sec:one-cell}
To better understand how the model problem we considered works and what does the ANN learn, we study a specific subset of the model problems defined in the previous section. Namely we consider a diffusion $\mu$ which is everywhere constant and equal to one apart from one cell of the mesh, where it is equal to $\mu_{MAX}=10^{\varepsilon}$ with $\varepsilon=2,4,8$. We stress that this family of problems was not present in the dataset we employed to train our ANN. Moreover, we define $d$ to be the distance between the cell with the largest diffusion coefficient $\mu_{MAX}$ from the center of the domain. Figure~\ref{fig:spectrum} (\textit{top-left}) shows a graphical representation of the problem. First we analyze the conditioning of the matrix $\mathrm A$. As expected it scales with respect to the mesh size $h$ as $h^{-2}$ and it is linearly proportional to the ratio between the maximum and minimum value of diffusion coefficient $\mu$ in the domain (i.e.\ is proportional to $\mu_{MAX}$). The position of the cell does not affect the condition number. There is just a small reduction of it when the cell touches the boundaries: this is simply due to the Dirichlet boundary conditions, see Figure~\ref{fig:spectrum} (\textit{bottom}). We then test our algorithm. Figure~\ref{fig:spectrum} (\textit{top-right}) shows that the optimal value of $\theta$ is almost always in the same range: between 0.5 and 0.9. Moreover, we see that the ANN is able to capture the overall relation between the normalized computational cost $t$ and the strong threshold. Hence, in all this cases, the predicted value $\theta^*$ by the AMG-ANN algorithm is near the true optimum. Finally, to gain some insight on the ANN we computed the feature maps extracted using the convolutional filters learned by the ANN. In Appendix~\ref{apx:cnn-feature-mapsl} we show some of them: we can see that from a layer to another the CNN is enhancing the features that deems relevant.

\begin{figure}[t]
        \centering
        \includegraphics[width=0.4\linewidth]{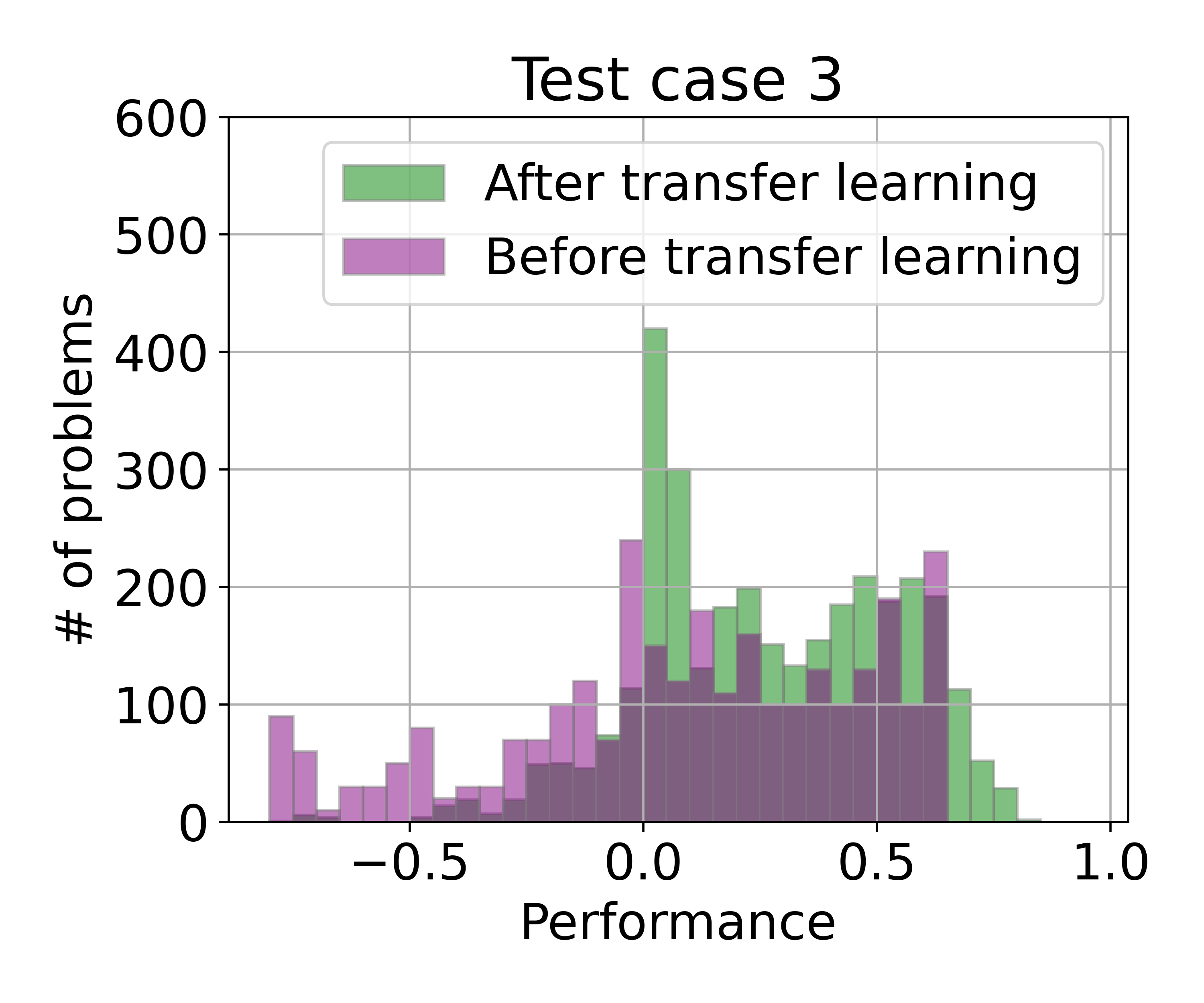}
        \caption{Relative performance improvement of the AMG method when applying our algorithm on Test case 3 before and after using transfer learning. The baseline is the ANN trained in Test case 2.}
        \label{fig:transfer-learn-hist}
\end{figure}

\subsection{Test Case 3: Linear Elasticity}
We consider the following linear elasticity problem
\begin{equation}
\begin{cases}
-\text{div}(\mathbb{C} \nabla_{S} \mathbf u)  = \mathbf f, \quad &\text{in} \: \Omega,\\
\mathbf u = \mathbf 0, &\text{on} \: \partial \Omega,
\end{cases}
\label{eq:elasticity}
\end{equation}
where $\Omega = (-1, 1)^3$ and $\nabla_S \mathbf u = \frac{1}{2} (\nabla \mathbf u + \nabla \mathbf{u}^\top)$ is the symmetric gradient and $\mathbb C$ is a rank-4 tensor that encodes the stress-strain relationship. Under the assumptions that the material under consideration is isotropic, 
 by introducing  the two Lam\'e coefficients $\lambda \in L^\infty(\Omega)$ and $\mu \in L^\infty(\Omega)$, the coefficient tensor reduces to
 \begin{linenomath}\begin{equation*}
\mathbb{C} \nabla_S \mathbf u = 2 \mu \nabla_S \mathbf u + \lambda \mathrm{tr}(\nabla_S \mathbf u) \mathrm{I}
\end{equation*}\end{linenomath}
where $\mathrm{tr}$ is the trace operator and $\mathrm{I}$ is the identity matrix. 
We can rewrite the Lam\'e parameters in terms of the Young's modulus $E > 0$
and the Poisson ratio $\nu \in (0, \frac{1}{2})$; we have $G = E/(1+\nu)$ and $\beta = \nu/(1-2\nu)$. The problem has been discretized by means of FE of order $p=1,2,3$ on a structured cartesian grid of diameter $h=\frac{\sqrt[3]{2}}{2}, ..., \frac{\sqrt[3]{2}}{2^{4-p}}$. We fix the Poisson ration $\nu = 0.29$ and choose a highly heterogeneous Young modulus $E$. Namely $E$ has the same pattern of the diffusion coefficient $\mu$ described in Section~\ref{sec:tc1} and is such that \texttt{size} is even. We also employ four different DoFs numbering to enhance the dataset. We solve each problem with 37 values of $\theta$ equispaced between 0.05 and 0.95. In this way, we build a dataset counting 5873 different problems and thus containing 217301 samples.

As we have done for the previous test case, we start with an architecture found like described in Section~\ref{sec:hyperparameters-study} and then we fine tuned the number of filters, the kernel size of the CNN and the width and depth of the dense part. The ANN has a (test) loss of $5.64 \cdot 10^{-4}$, Figure~\ref{fig:loss-history} shows the history of the training loss. Evaluated on the test dataset, the algorithm has an accuracy of 90.71\% and an average reduction of the computational cost of almost 23\% in average. Details about its performance and its scaling are shown in Figure~\ref{fig:perf-and-cost} and Table~\ref{tab:perf_summary}. The number of iterations needed reach convergence is between 9 and 148 with $\theta=0.5$ and between 8 and 87 with $\theta=\theta^*$.  There is still margin of improvement concerning the reduction of the cost since $P_{MAX}$ is relatively small, however it is possible to see a remarkable reduction of the scaling of the computational cost in the bottom row of Figure~\ref{fig:perf-and-cost}. 

\paragraph{Speeding up the training} We show that transfer learning ad layer freezing can be used to to speed up the training phase. We use the ANN we trained for Test case 2 as a pretrained network. Namely, we use it as the starting point of the training but we keep the weights of the convolutional block frozen. In this way we significantly reduce the cost of training the network since many parameters do not need the computation of the gradient. By caching the representation learnt by the CNN we can lower even more the cost of the training. Before the transfer learning we have $PB=59.32\%$ and $P_M=6.22\%$, which is slightly better than using $\theta=0.5$. After the training, the loss is $3.97 \cdot 10^{-4}$, details about the performance are shown in Figure~\ref{fig:transfer-learn-hist} and Table~\ref{tab:perf_summary}. With a computational cost that is more than halved we have trained an ANN that performs better than the one trained from scratch.

\begin{remark}

A key aspect of the proposed AMG-ANN algorithm is how much data and variation in the data are required in order to obtain a good estimate of $\theta^*$ for your problem. Unfortunately, this is almost impossible to determine a-priori, since it is equivalent to determine how much an ANN generalizes, which is still an open problem of the field of DL. However, from the presented test cases we can conclude the following: 
\begin{itemize}
    \item The algorithm starts to struggle when making predictions on geometries $\Omega$ that has never seen before. It still proves effective if $\Omega$ is the union of geometries present in the dataset (see Section~\ref{sec:unseen-domains} for comparing the performance on ``Holes'' which is a replica of ``Plate with hole'' ).
    \item The algorithm does not perform well when applied to problems with a different physics from the one in the training dataset (e.g., making prediction on linear elasticity problem using the dataset of Test case 2). However, transfer learning reduces significantly the training cost on new data. 
\end{itemize}

\end{remark}

\section{Conclusions} \label{sec:conclusions}
We extended the algorithm presented in \cite{antonietti2023accelerating} by adding features that allow it to encompass the complexity of three-dimensional problems. We discussed a pre-processing stage to ensure data quality, an augmented the architecture of the ANN and an enhanced version of the pooling for sparse matrices. Additionally, we demonstrated a method for accelerating training by implementing layer freezing techniques.
This allows the algorithm to perform well also on a wide class of challenging problems coming from the FE discretization of elliptic PDEs with rough coefficients. We have provided numerical results for three different test cases concerning a Poisson problem with a highly heterogeneous diffusion coefficient on (\textit{i}) structured and (\textit{ii}) unstructured grids and (\textit{iii}) a linear elasticity problem with a highly heterogeneous Young's modulus. Our algorithm performs better, other than using the standard literature value, in more than 90\% of the cases and reduces significantly the computational cost (up to 30\% on average). This is enabled by the fact that after the expensive training phase (which is done just one), the online cost of the algorithm is negligible.

The algorithm hinges upon the pooling operator to compress a large sparse matrix into a small tensor to feed into an ANN. Indeed, we have proved that the pooling operator reduces the computational cost and preserves relevant information needed to perform the complex regression task at hand. Moreover, our algorithm can be introduced with minimal changes into any existing code that uses an AMG solver.

In future works we plan to generalize our algorithm by using one single ANN to make predictions on linear systems stemming from different underlying PDEs and being able to predict the optimal value of other parameters of the AMG method such as the type of smoother and the number of smoothing steps. We also plan to try more advanced computer vision models, such as transformers \cite{dosovitskiy2020image}. Finally, we aim to further investigate the properties of the pooling operator applied to sparse matrices and gain a better understanding of what does the CNN learn.

\paragraph*{Funding}
{\small
M.C., P.F.A and L.D. are members of the INdAM Research group GNCS. P.F.A has been partially funded by the research projects PRIN17 (n. 201744KLJL) and PRIN 2020 (n. 20204LN5N5) funded by Italian Ministry of University and Research (MUR) and partially funded by ICSC—Centro Nazionale di Ricerca in High Performance Computing, Big Data, and Quantum Computing funded by European Union—NextGenerationEU. L.D. has been partially funded by the research project PRIN 2020 (n. 20204LN5N5) funded by MUR. L.D. acknowledges the support by the FAIR (Future Artificial Intelligence Research) project, funded by the NextGenerationEU program within the PNRR-PE-AI scheme (M4C2, investment 1.3, line on Artificial Intelligence), Italy. The present research is part of the activities of ``Dipartimento di Eccellenza 2023-2027''.}

\paragraph*{CRediT authorship contribution statement}
{\small P.F.A. and L.D.: Conceptualization, Methodology, Review and editing, Supervision, Project administration, Funding acquisition. M.C.: Conceptualization, Methodology, Software, Validation, Formal analysis, Investigation, Data curation, Visualization, Original draft.}
\paragraph*{Declaration of competing interests} 
{\small The authors declare that they have no known competing financial interests or personal relationships that could have appeared to influence the work reported in this paper.}

\appendix
\section{Basic concepts of Deep Learning}\label{apx:dl}
A deep learning regression model is a function $\mathscr{F} \, :\, \mathbb{R}^{N_{IN}} \rightarrow \mathbb{R}^{N_{OUT}}$ that maps an input $\mathbf{x}$ to an output $\tilde{\mathbf{y}}$ and depends on a vector of parameters ${\boldsymbol{\gamma}}$. The parameters ${\boldsymbol{\gamma}}$ are chosen by minimizing the error (loss) evaluated on a training dataset of known input-output couples $\{(\mathbf{x}^i,\mathbf{y}^i)\}_{i=1}^{T}$, namely
\begin{linenomath}\begin{equation*}
\min_{\boldsymbol{\gamma}} \frac{1}{T}\sum_{i=1}^T \mathscr{L}(\mathbf y^i, \mathscr F(\mathbf x^i, \boldsymbol{\gamma}))
\end{equation*}\end{linenomath}
where $\mathscr{L}$ is the loss function. Namely, we employ the mean squared error (MSE). The optimization is performed by means of gradient descent updates that uses a mini-batch of the training dataset to compute $\nabla_{\boldsymbol{\gamma}} \mathscr{F}$ with automatic differentiation.

\paragraph{Feed-forward neural networks}
Let $\mathbf y^{(0)} = \mathbf x$ and $\mathbf y^{(L)} = \tilde{\mathbf{y}}$, a dense feed-forward neural network (FFNN) of depth $L$ is the composition of $L$ functions called layers defined as

\begin{linenomath}\begin{equation*}
\mathbf y^{(l)} = h^{(l)}(\mathrm W^{(l)} \mathbf y^{(l-1)} + \mathbf b^{(l)}), \quad l=1,...,L
\end{equation*}\end{linenomath}
where $\mathrm{W}^{(l)} \in \mathbb{R}^{N_l \times N_{l-1}}$ (weights) and $\mathbf{b}^{(l)} \in \mathbb{R}^{N_l}$ (biases) are the parameters $\boldsymbol{\gamma}$, and $h^{(l)}$ is a scalar non-linear activation function applied component-wise.

\paragraph{Convolutional neural networks}
Convolutional neural networks (CNNs) are neural networks that use convolution instead of matrix
multiplication in at least one of their layers. They have great success in computer vision tasks thank to their ability to exploit the structured data format of images. Indeed, convolution layers enjoy three properties: shared parameters, that is, each parameter is tied to multiple component of the input; sparse interactions, that is, each component of $\mathbf y^{(l)}$ depends only on a subset of the components of $\mathbf y^{(l -1)}$, (this also entails a lower number of parameters and thus a greater efficiency); equivariance
to translation, that is the application of a translation and a convolution can be interchanged to obtain the same result. The last ingredient of a CNNs layer is the pooling operator. After several parallel convolutions and a non-linear activation, we replace the output of the layer at a certain location with a
summary statistic of the nearby outputs. A popular option is the max pooling \cite{krizhevsky2017imagenet}, which reports the maximum value over a rectangular neighborhood. The aim of the pooling operation is to reduce the computational cost by limiting the number of parameters in subsequent layers, learn hierarchical representations of features, where lower layers capture low-level details, and higher layers capture more abstract and complex features, and adding invariance to small translations.

\section{Pooling visualization}\label{apx:pooling_visual}
\begin{figure}[H]
    \makebox[\textwidth][c]{\includegraphics[trim={1.5cm 0.8cm 0.5cm 0.7cm},clip,width=0.95\textwidth]{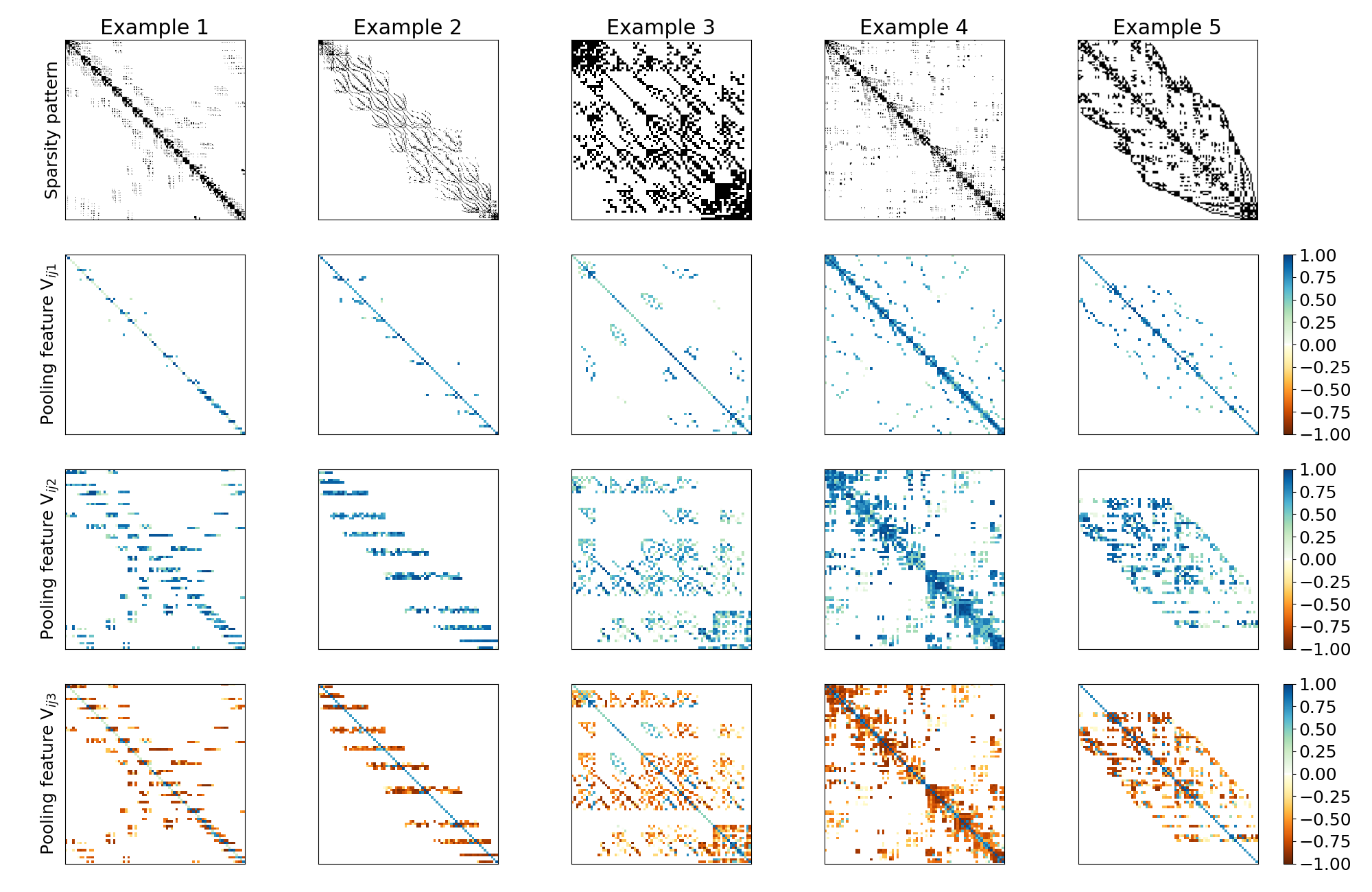}}
    \caption{Visual representation of the normalized pooling results $\hat{\mathrm V}$ that is fed into the ANN for some exemplary problems taken from the dataset described in Section~\ref{sec:tc2}.}
    \label{fig:views}
\end{figure}

\section{CNN feature maps visualization}\label{apx:cnn-feature-mapsl}
\begin{figure}[H]
\setlength{\tabcolsep}{0pt}
    \centering
        \begin{tabular}{cccc}
             \includegraphics[width=0.245\linewidth]{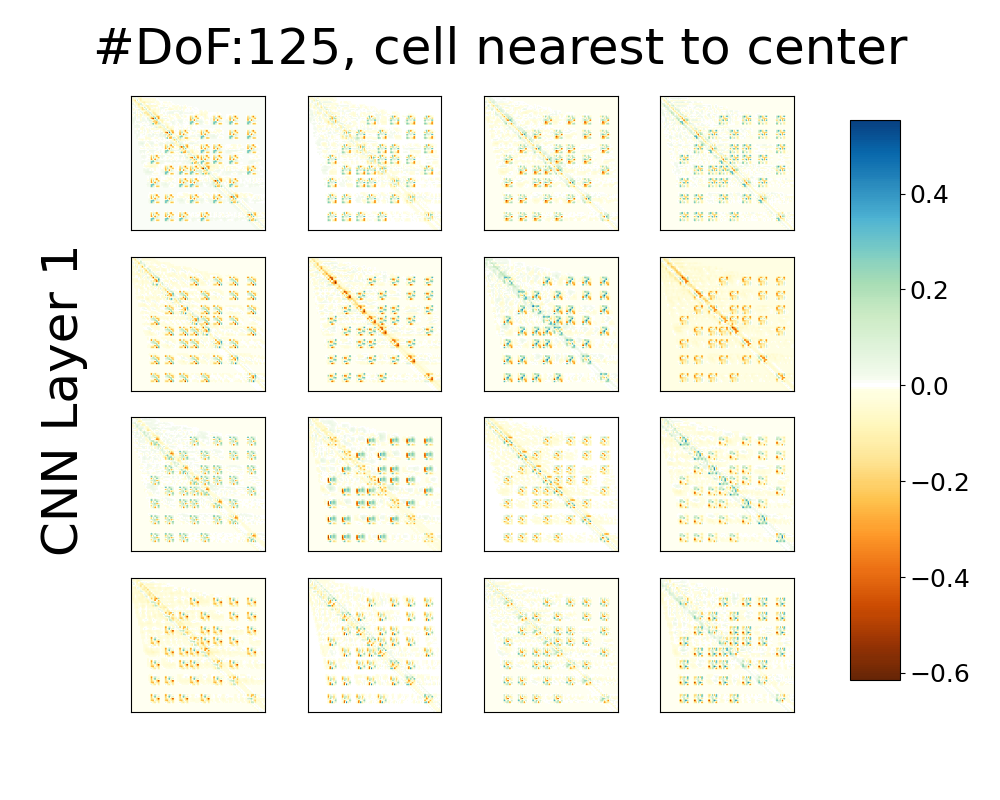}
            & \includegraphics[width=0.245\linewidth]{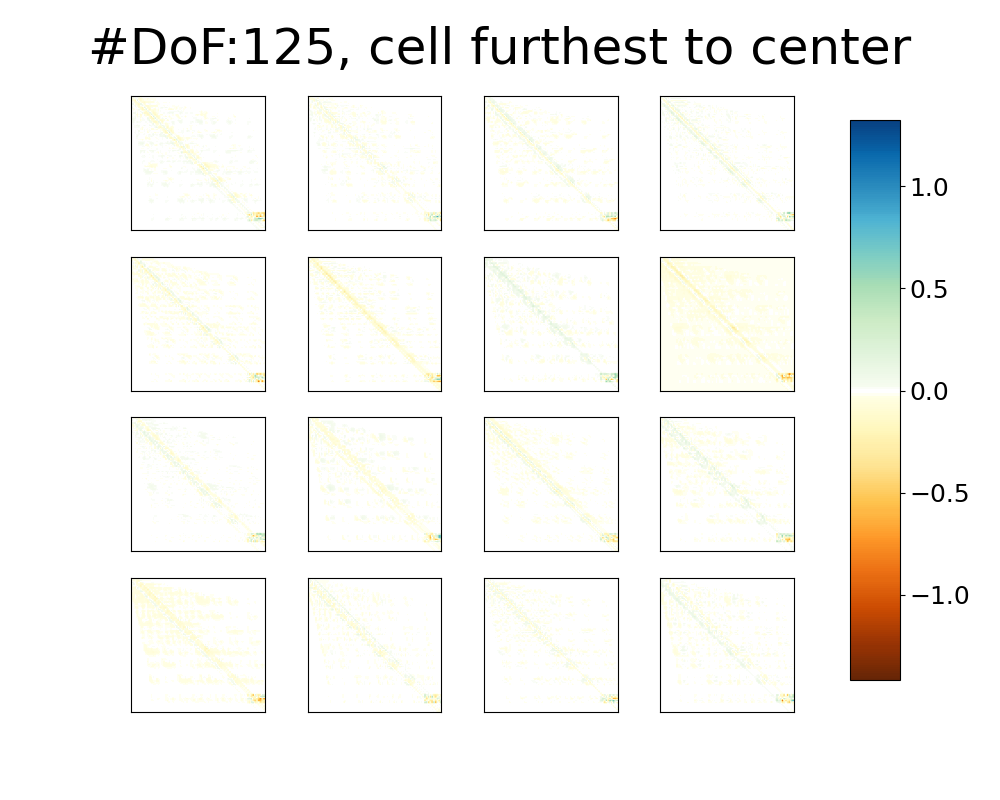}
            & \includegraphics[width=0.245\linewidth]{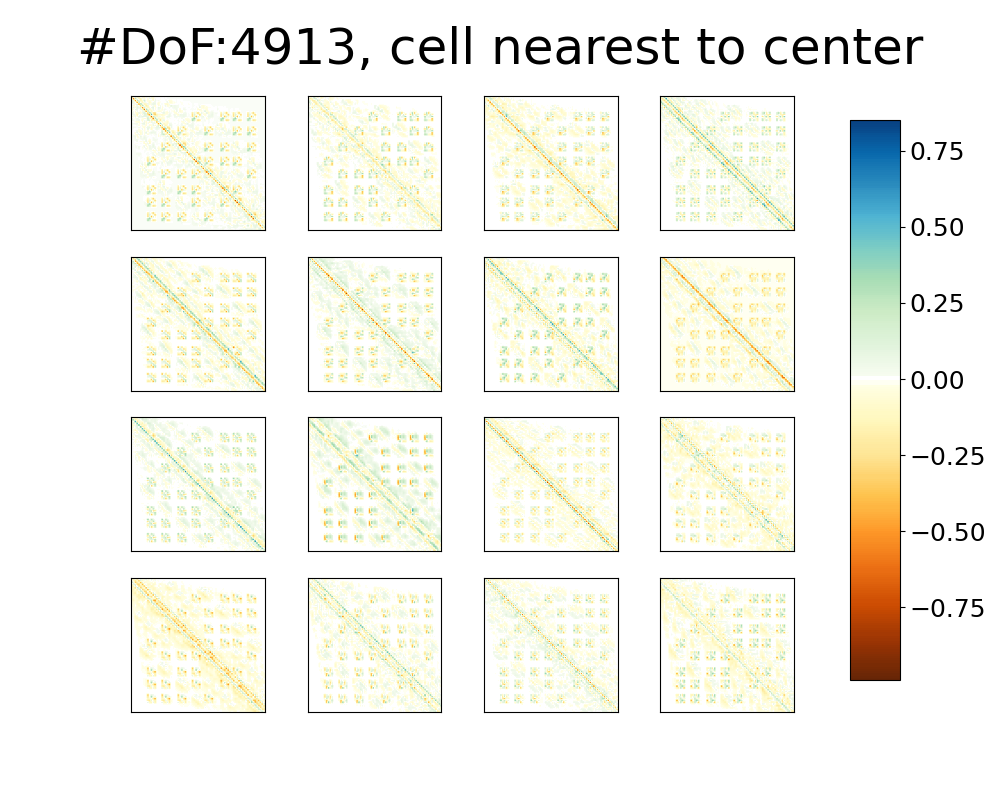}
            & \includegraphics[width=0.245\linewidth]{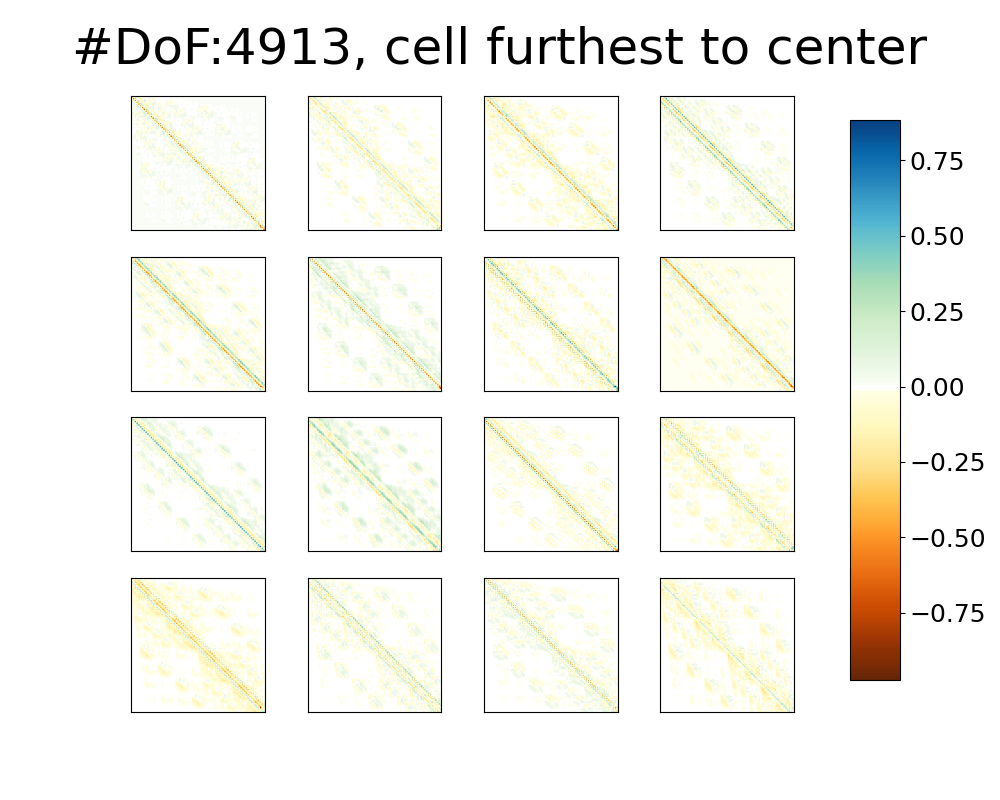}\\[-4pt]
            
             \includegraphics[width=0.245\linewidth]{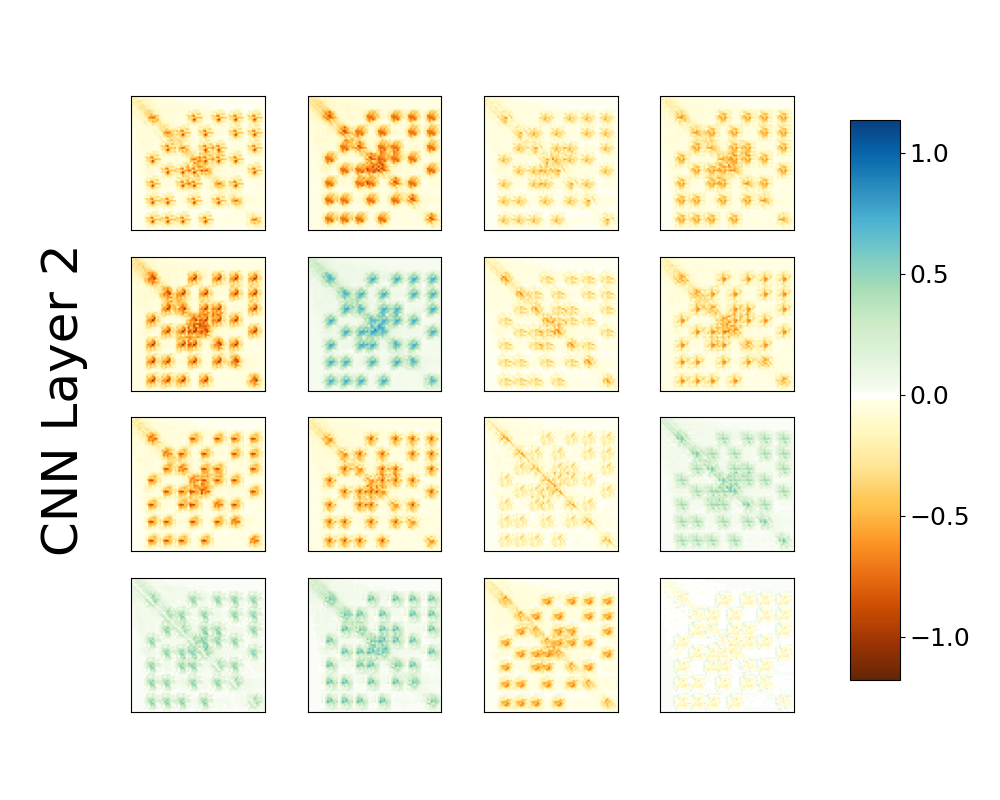}
            & \includegraphics[width=0.245\linewidth]{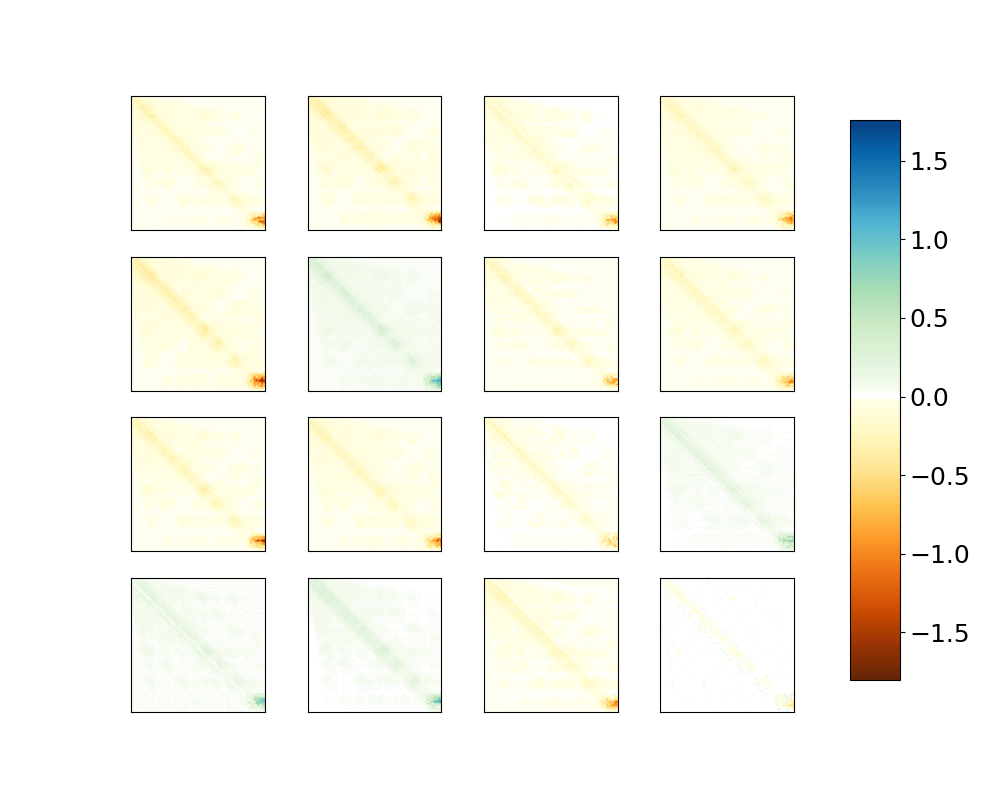}
            & \includegraphics[width=0.245\linewidth]{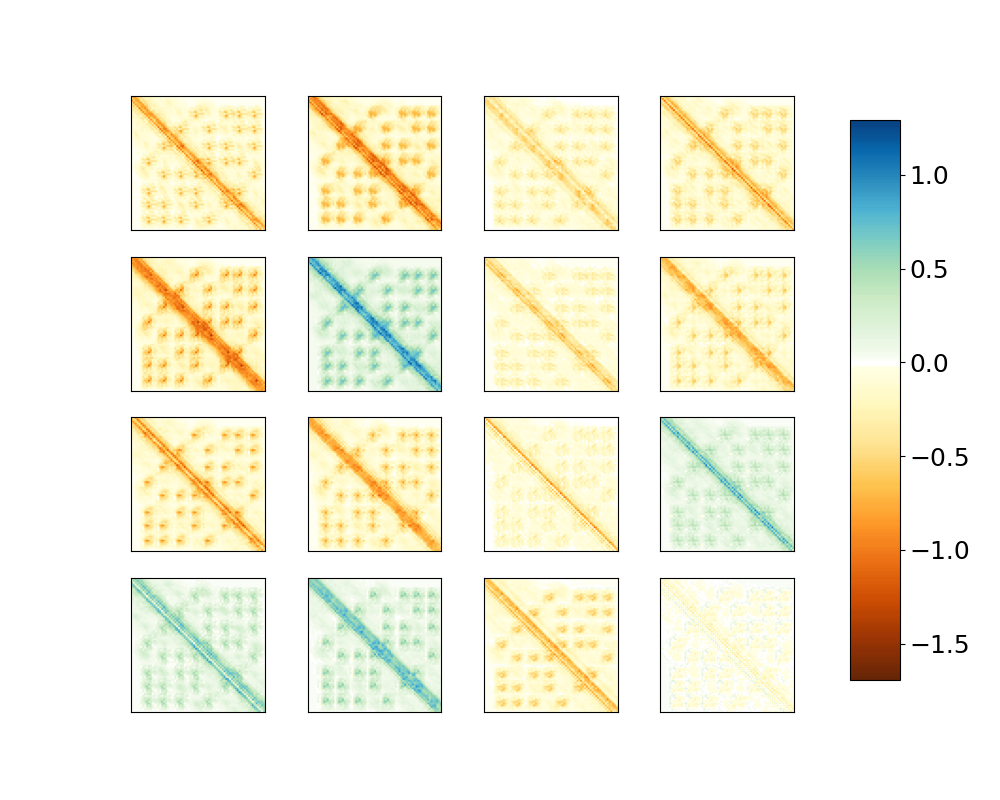}
            & \includegraphics[width=0.245\linewidth]{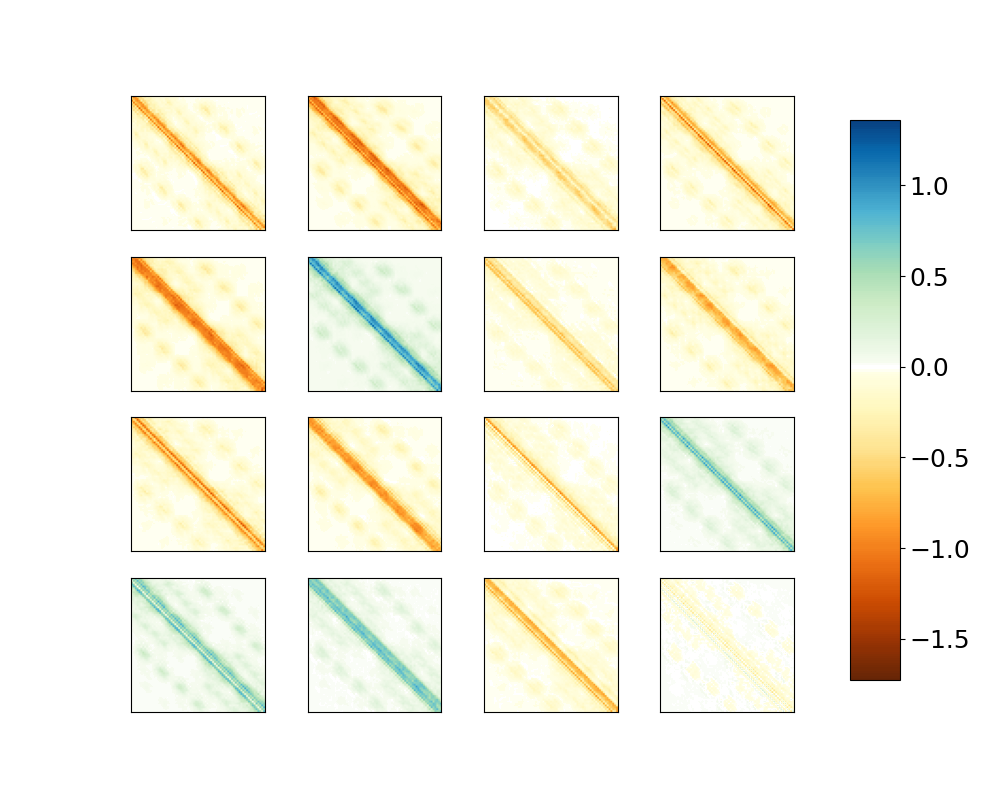}\\[-4pt]

             \includegraphics[width=0.245\linewidth]{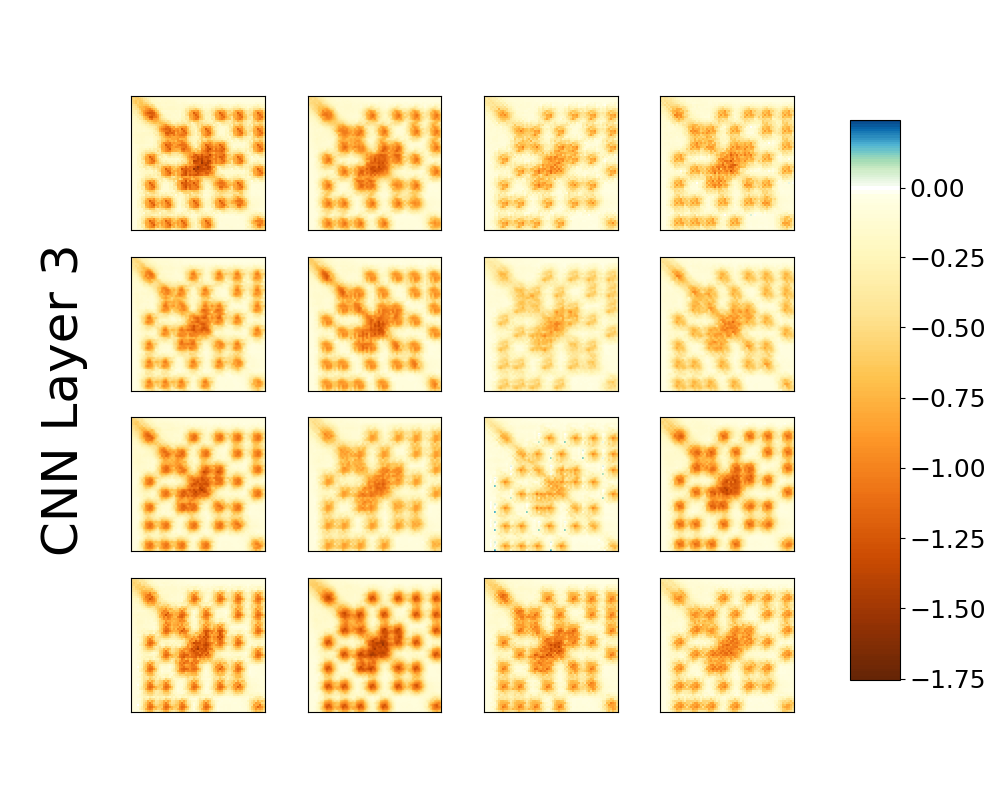}
            & \includegraphics[width=0.245\linewidth]{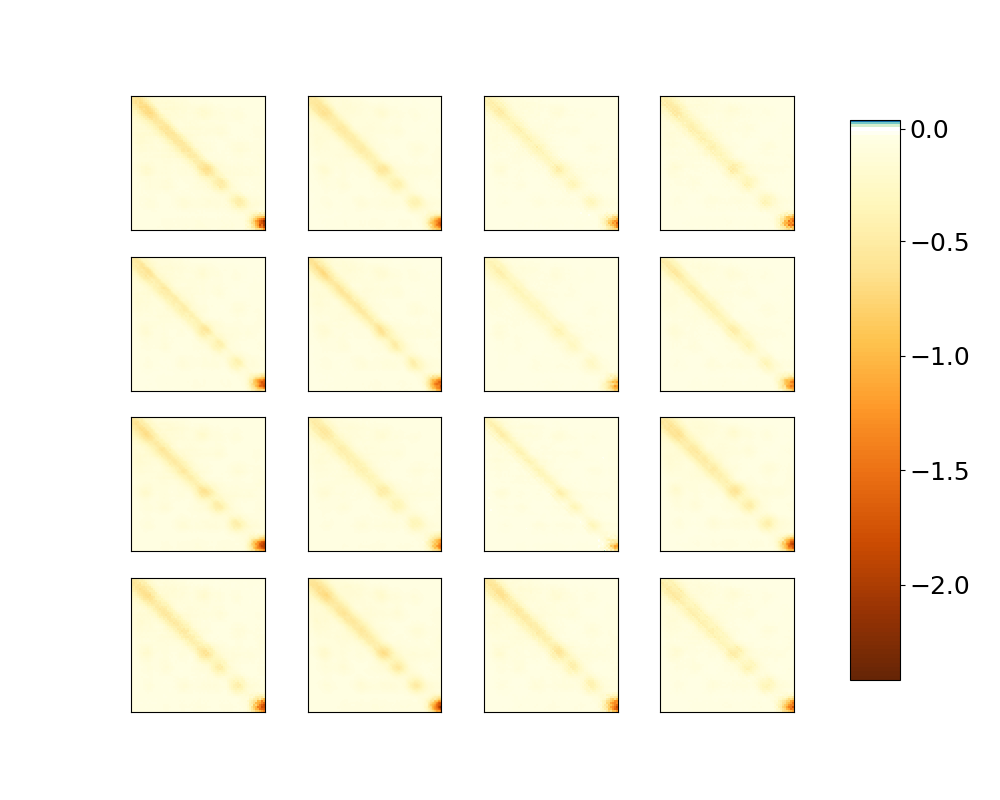}
            & \includegraphics[width=0.245\linewidth]{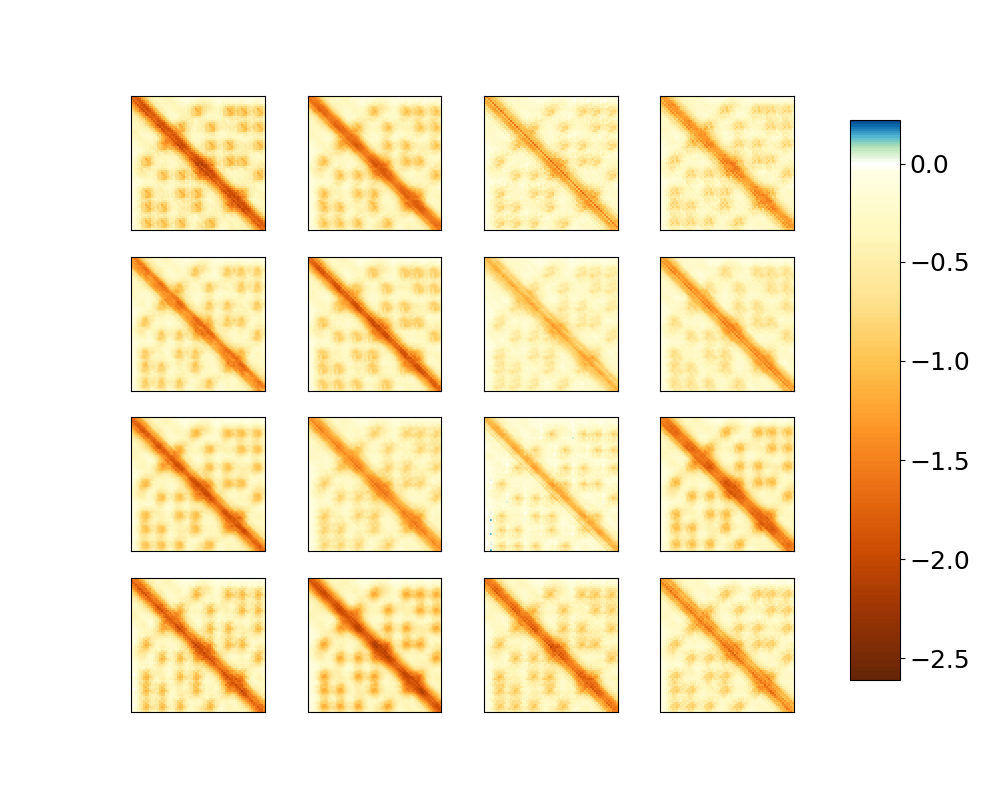}
            & \includegraphics[width=0.245\linewidth]{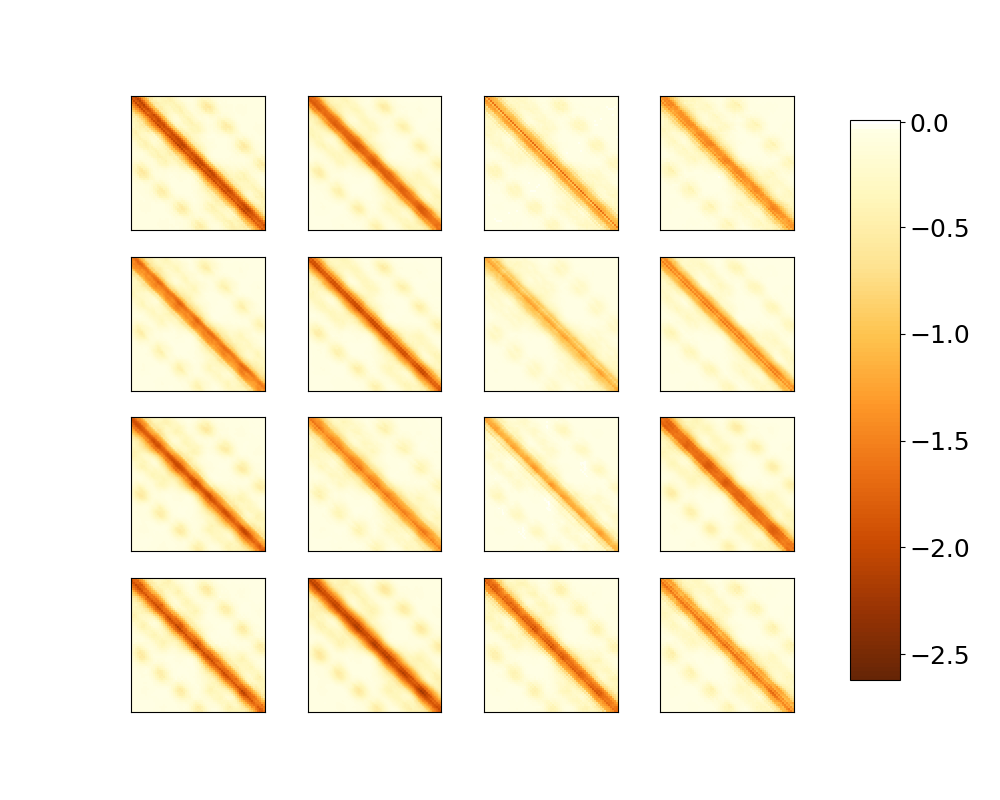}\\[-4pt]
        \end{tabular}%
    \caption{Example of feature maps of the CNN layers for the problems described in Section~\ref{sec:one-cell}.\label{fig:cnn-feature-maps}}
\end{figure}
\bibliographystyle{abbrv}
\bibliography{references}

\end{document}